\documentclass[12pt,draftcls,onecolumn]{IEEEtran}

\usepackage{amssymb}
\usepackage{makeidx}         % allows index generation
\usepackage{amsmath}
\usepackage{multicol}        % used for the two-column index
\usepackage[bottom]{footmisc}% places footnotes at page bottom
\usepackage{graphicx}
\usepackage{color}
\usepackage{amsfonts}
\usepackage[latin1]{inputenc}
\usepackage{epstopdf}
\usepackage{url}
\usepackage[usenames,dvipsnames]{xcolor}

\newcommand{\eq}{\begin{equation}}
\newcommand{\eeq}{\end{equation}}
\newcommand{\eqn}{\begin{eqnarray}}
\newcommand{\eeqn}{\end{eqnarray}}

\newcommand{\bsea}{\begin{subeqnarray}}
\newcommand{\esea}{\end{subeqnarray}}
\newcommand{\nn}{\nonumber}

\newcommand{\Sp}[2]{\left< #1,#2 \right> }

\newcommand{\tr}{\mathop{\rm tr}}  %traccia

\newcommand{\Ac}{ \mathcal{A}}

\newcommand{\Dc}{ \mathcal{D}}

\newcommand{\Gc}{ \mathcal{G}}

\newcommand{\Lc}{ \mathcal{L}}

\newcommand{\Nc}{ \mathcal{N}}

\newcommand{\Qc}{ \mathcal{Q}}

\newcommand{\Sc}{ \mathcal{S}}

\newcommand{\Vc}{ \mathcal{V}}

\newcommand{\Xc}{ \mathcal{X}}

\newcommand{\Cs}{ \mathbb{C}}
\newcommand{\Ds}{ \mathbb{D}}
\newcommand{\Es}{ \mathbb{E}}

\newcommand{\Ns}{ \mathbb{N}}

\newcommand{\Rs}{ \mathbb{R}}

\newcommand{\Zs}{ \mathbb{Z}}

\def\qed{\hfill \vrule height 7pt width 7pt depth 0pt \smallskip}

\newcounter{pippo}

\newtheorem{remark}{Remark}[section]
\newtheorem{teor}{Theorem}[section]
\newtheorem{corr}{Corollary}[section]
\newtheorem{propo}{Proposition}[section]
\newtheorem{lemm}{Lemma}[section]
\newtheorem{exam}{Example}
\newtheorem{probl}[pippo]{Problem}
\newtheorem{defn}{Definition}[section]
\newcommand{\proof}{\noindent {\bf Proof. }}
\newcommand{\teo}{\begin{teor}}
\newcommand{\eteo}{\end{teor}}
\newcommand{\cor}{\begin{corr}}
\newcommand{\ecor}{\end{corr}}
\newcommand{\prop}{\begin{propo}}
\newcommand{\eprop}{\end{propo}}
\newcommand{\lem}{\begin{lemm}}
\newcommand{\elem}{\end{lemm}}
\newcommand{\ex}{\begin{exam}}
\newcommand{\eex}{\end{exam}}
\newcommand{\pb}{\begin{probl}}
\newcommand{\epb}{\end{probl}}
\newcommand{\df}{\begin{defn}}
\newcommand{\edf}{\end{defn}}
\newcommand{\aprop}{\begin{apropo}}
\newcommand{\eaprop}{\end{apropo}}
\newcommand{\alem}{\begin{alemm}}
\newcommand{\ealem}{\end{alemm}}
\newcommand{\rem}{\begin{remark}}
\newcommand{\erem}{\end{remark}}

\newcommand{\Stp}  {\mathcal{S}_m}

\newcommand{\Sph}  {\mathcal{S}_{m+l}}
\newcommand{\botc}  {\;\perp\;}

\newcommand{\Dd}  {\mathrm{D}}
\newcommand{\Ld}  {\mathrm{L}}
\newcommand{\Pd}  {\mathrm{P}}
\newcommand{\Td}  {\mathrm{T}}

\newcommand{\Mb}  {\mathbf{M}}
\newcommand{\Qb}  {\mathbf{Q}}

\begin{document}

{\color{black}

\title{AR Identification of\\ Latent-variable Graphical Models}

\author{Mattia~Zorzi, Rodolphe~Sepulchre\thanks{This paper presents research
results of the Belgian Network DYSCO (Dynamical Systems, Control, and
Optimization), funded by the Interuniversity Attraction Poles Programme, initiated
by the Belgian State, Science Policy Office. The scientific responsibility
rests with its authors. This research is also supported by FNRS (Belgian Fund for Scientific Research).} \thanks{M. Zorzi is with the
Department of Electrical Engineering and Computer Science,
University of Liège, 4000 Liège, BE,
({\tt\small mzorzi@ulg.ac.be})} \thanks{R. Sepulchre is with the
Department of Engineering, University of Cambridge, Cambridge CB2 1PZ, UK,
({\tt\small r.sepulchre@eng.cam.ac.uk}), and the
Department of Electrical Engineering and Computer Science,
University of Liège, 4000 Liège, BE.} }

\markboth{DRAFT}{Shell \MakeLowercase{\textit{et al.}}: Bare Demo of IEEEtran.cls for Journals}

\maketitle

\begin{abstract}
The paper proposes an identification procedure for autoregressive gaussian stationary stochastic processes wherein the manifest (or observed) variables are mostly
related through a limited number of latent (or hidden) variables. The method  exploits the sparse plus low-rank decomposition of the inverse of the manifest spectral density and the efficient convex relaxations recently
 proposed for such decomposition.
\end{abstract}


\section{Introduction}

Gaussian processes and their representation by graphical models have gained popularity through science and engineering, \cite{LAURITZEN_1996,Willsky02multiresolutionmarkov}. The objective of the present paper is to derive an identification procedure for gaussian stochastic processes  whose manifest (observed) variables are correlated primarily through a restricted number of latent (hidden) variables. Here, (the few) latent variables are fictitious elements introduced by the modeler. The resulting graphical model (or equivalently latent-variable  graphical model) has a two layer structure, one layer for the manifest (observed) nodes and one layer for the latent (hidden) nodes. The hope is that in many applications of interest, the few extra nodes  in the hidden layer allow for a drastic  reduction of edges in the observed layer, because the observed nodes become nearly independent when conditioned to the hidden nodes. As a consequence, allowing for latent variables in the identification of the stochastic model may improve scalability and robustness of the algorithm. This paradigm was exploited in the framework of gaussian random vectors in the recent paper \cite{Chandrasekaran_latentvariable}. The authors exploited the sparse plus low-rank (S+L) decomposition of the manifest concentration matrix (the inverse of the covariance matrix corresponding to the manifest variables) to propose an efficient formulation of the identification problem.

The present paper focuses on the generalization of this approach to autoregressive (AR) gaussian stationary processes, exploiting the analog sparse plus low-rank decomposition of the inverse of the manifest spectral density
(the spectral density of the manifest variables). It thereby connects the extensive recent research on convex regularization of sparsity and low-rank constraints \cite{Chandrasekaran_latentvariable,CHANDRASEKARAN_RANK_sparsaity_2011,S_PLUS_L_OTAZO,Choi2010GAUSSIANMULTIRES,SAUNDERSON_TREE_STRUCTRURED}
to the classical covariance extension approach for the identification of gaussian stationary processes \cite{burg1975maximum,A_NEW_APPROACH_BYRNES_2000}. It also provides a generalization of recent contributions that introduced sparsity constraints (but no latent variables) in the identification of autoregressive processes \cite{SONGSIRI_TOP_SEL_2010,SONGSIRI_GRAPH_MODEL_2010,ARMA_GRAPH_AVVENTI}.

The paper is organized as follows. After mathematical preliminaries, Section  \ref{section_pb_formulation}
introduces the main ideas of the proposed identification scheme in non technical terms. The identification of the graphical model and the identification of the autoregressive model are formulated as two distinct optimization problems. The first one uses sparsity and low-rank regularizers to recover the model structure. It is further analyzed in Section \ref{sec_top_sel}. The second  one solves an exact covariance extension problem  for a fixed graphical model. It is further  analyzed in Section \ref{sec_cov_ext}. Finally, in Section \ref{sec_id_procedure}
we discuss an illustrative example and test our method to international stock return data.

\subsection*{Notation}
We endow the vector space $\Rs^{m\times m}$ with the usual inner product $\Sp{X}{L}=\tr(X L^T)$. $\Qb_{m}$ denotes the vector space of symmetric matrices of dimension $m$, if $X\in\Qb_m$ is
positive definite (semi-definite) we write $X\succ 0$ ($X\succeq 0$). A matrix $A\in\Rs^{l\times m(n+1)}$ with $l\leq m$ will be partitioned as
$A=\left[
     \begin{array}{cccc}
      A_0  & A_1 & \ldots & A_n \\
     \end{array}
   \right] $ with $A_j\in\Rs^{l\times m}$. $\Mb_{m,n}$ is the vector space of matrices $Y:=\left[
                                                                                                         \begin{array}{cccc}
                                                                                                          Y_0 & Y_1 & \ldots & Y_{n} \\
                                                                                                         \end{array}
                                                                                                       \right]$ with  $Y_0\in\Qb_m$ and $Y_1\ldots Y_{n}\in\Rs^{m\times m}$. The corresponding inner product is $\Sp{Y}{Z}=\tr(YZ^T)$. The linear mapping $\Td: \Mb_{m,n} \rightarrow \Qb_{m(n+1)}$ constructs a symmetric {\em Toeplitz}
matrix from its first block row in the following way:\eq \label{definizione_operatore_T}\Td(Y)=\left[
                    \begin{array}{cccc}
                      Y_{0} & Y_{1} & \ldots & Y_{n} \\
                      Y_{1}^T & Y_{0} & \ddots & \vdots \\
                      \vdots & \ddots &  \ddots& Y_1\\
                      Y_{n}^T & \ldots & Y_1^T & Y_{0} \\
                    \end{array}
                  \right].
\eeq The adjoint of $\Td$ is a mapping $\Dd: \Qb_{m(n+1)}\rightarrow \Mb_{m,n}$ defined as follows. If $X \in\Qb_{m(n+1)}$ is partitioned as
\eq \label{partizione_X}X=\left[
                    \begin{array}{cccc}
                      X_{0 0} & X_{0 1} &\ldots & X_{0n} \\
                      X_{01}^T & X_{11} &\ldots & X_{1n} \\
                        \vdots  & \vdots &   & \vdots \\
                      X_{0 n}^T & X_{1 n}^T & \ldots  & X_{nn} \\
                    \end{array}
                  \right]\eeq
then $\Dd(X)=\left[
                      \begin{array}{ccc}
                        \Dd_0(X) & \ldots & \Dd_{n}(X) \\
                      \end{array}
                    \right]
$ where \eq \label{operatoreD}\Dd_0(X)=\sum_{h=0}^{n} X_{h h},\,\; \Dd_j(X)=2\sum_{h=0}^{n-j} X_{h\; h+j},\; \; j=1\ldots n.\eeq
We define the index set $E_m \subseteq V_m\times V_m$ with $V_m:=\{1,2,\ldots m\}$, and its complement set is denoted by $E_m^c$. The cardinality of $E_m$ is denoted by $|E_m|$.
The projection map $\Pd_{E_m}: \Rs^{m\times m} \rightarrow \Rs^{m\times m} $ is defined as follows
\eq \label{operatore_P_E}\Pd_{E_m}(X)=\left\{
                  \begin{array}{ll}
                    (X)_{kh}, & (k,h)\in E_m \\
                    0, & \hbox{otherwise}
                  \end{array}
                \right.
\eeq where $(X)_{kh}$ is the entry of $X$ in position $(k,h)$. Similarly, $\Pd_{E_m}(Y)$ with $Y\in\Mb_{m,n}$ denotes \eq \label{def_op_T_blocchi}\left[
                                                                        \begin{array}{cccc}
                                                                          \Pd_{E_m}(Y_0) & \Pd_{E_m}(Y_1) & \ldots & \Pd_{E_m}(Y_n) \\
                                                                        \end{array}
                                                                      \right].\eeq

Functions on the unit circle $\{e^{i\vartheta} \hbox{ s.t. } \vartheta \in[-\pi,\pi]\}$ will be denoted by capital Greek letters, e.g. $\Phi(e^{i\vartheta})$ with $\vartheta\in[-\pi,\pi]$, and the dependence upon $\vartheta$ will be dropped if not needed, e.g. $\Phi$ instead of $\Phi(e^{i\vartheta})$. $\Ld_2^{m\times m}$ denotes the space of $\Cs^{m\times m}$-valued functions defined on the unit circle which are square integrable. Given $\Phi\in\Ld^{m\times m}_2$, the shorthand notation $\int \Phi$ denotes the integration of $\Phi$ taking place on the unit circle with respect to the normalized {\em Lebesgue} measure. Then, the inner product in $\Ld_2^{m\times m}$ is $\Sp{\Phi}{\Sigma}=\tr\int \Phi\Sigma^*$. Similarly, $\Pd_{E_m}: \Ld_2^{m\times m} \rightarrow \Ld_2^{m\times m}$ is defined as in (\ref{operatore_P_E}) where $X$ is replaced by $\Phi(e^{i\vartheta})$. Moreover, $\sigma_k(\Phi(e^{i\vartheta}))$
denotes the $k$-th largest singular value of $\Phi(e^{i\vartheta})$ at $\vartheta$, i.e. $\sigma_1(\Phi(e^{i\vartheta}))\geq \sigma_2(\Phi(e^{i\vartheta}))\geq \ldots \geq \sigma_m(\Phi(e^{i\vartheta}))$ for each $\vartheta\in [-\pi,\pi]$. $\Ac_m$ denotes the linear space of $\Cs^{m\times m}$-valued analytic functions on the unit circle. Given $\Lambda\in\Ac_m$, we define the
norm \eq \| \Lambda\|=\sup_{\vartheta\in[-\pi,\pi]} \sigma_1(\Lambda(e^{i\vartheta}))\eeq
and the (normal) rank \eq \mathrm{rank}(\Lambda):=\underset{\vartheta \in[-\pi,\pi]}{ \max}\;\mathrm{rank} (\Lambda(e^{i\vartheta})).\eeq
If $\Phi(e^{i\vartheta})$ is positive definite (semi-definite) for each $\vartheta\in[-\pi,\pi]$, we will write $\Phi\succ 0$ ($\Phi \succeq 0$).
$\Stp$ denotes the family of functions $\Phi$ such that $\Phi=\Phi^*$ and $c_1 I\preceq \Phi\preceq c_2 I$ for some $c_1,c_2>0$.
We define the following family of matrix pseudo-polynomials
\eqn \Qc_{m,n}&=&\left\{\sum_{j=-n}^n e^{-ij\vartheta}R_j \hbox{ s.t. } R_j=R_{-j}^T \in\Rs^{m\times m}\right\}.\eeqn
The {\em shift operator} is defined as \eq \Delta(e^{i\vartheta}):=\left[
                                                                          \begin{array}{cccc}
                                                                            I_m & e^{i\vartheta}I_m & \ldots & e^{in \vartheta} I_m \\
                                                                          \end{array}
                                                                        \right].
\eeq
Given $X\in\Qc_{m(n+1)}$, by direct computation we get
\eqn \label{new_reparametrization} && \hspace{-0.8cm}\Delta(e^{i\vartheta}) X\Delta(e^{i\vartheta})^*\nn\\ && =\Dd_0(X)+\frac{1}{2}\sum_{j=1}^n e^{-ij\vartheta}\Dd_j(X)+e^{ij\vartheta}\Dd_j(X)^T,\eeqn
therefore $\Delta X \Delta^*\in\Qc_{m,n}$. On the other hand, any element in $\Qc_{m,n}$ may be parameterized as (\ref{new_reparametrization}) because $\Dd$ is a surjective map. We conclude that
\eq \label{set_Qc_mn_reparametrized}\Qc_{m,n}=\{\Delta X \Delta^* \hbox{ s.t. } X\in\Qb_{m(n+1)}\}.\eeq

\section{Problem Formulation}\label{section_pb_formulation}

\subsection{AR Model Identification}
Let $L_2^m(\Omega, \Ac,P)$ be the Hilbert space of second order $\Rs^m$-valued gaussian random vectors defined in the probability space $\{\Omega, \Ac,P\}$.
An $\Rs^m$-valued gaussian stochastic process $x^m$ is an ordered collection of random vectors $x^m=\{x^m(t);\; t\in\Zs\}$ in $L_2^m(\Omega, \Ac,P)$. Moreover, we assume $x^m$ is zero mean, stationary and purely nondeterministic. It is completely described by its spectral density
\eq \Phi_m(e^{i\vartheta})=\sum_{j=-\infty}^\infty e^{-ij \vartheta} R_j\eeq
where $R_j:=\Es[x^m(t+j) x^m(t)^T]$ denotes the $j$-th covariance lag. An empirical estimate $\hat R_j$ of $R_j$ is computed from a finite-length realization of $x^m$, i.e. $\mathrm{x}^m(1),\mathrm{x}^m(2),\ldots \mathrm{x}^m(N)$, as follows
\eq \label{estimated_covlag}\hat R_j=\frac{1}{N}\sum_{t=0}^{N-j} \mathrm{x}^m(t+j)\mathrm{x}^m(t)^T.\eeq
The estimate $\hat \Phi_m^\circ$ of $\Phi_m$
that maximizes the {\em entropy rate}, \cite{COVER_THOMAS}, and that matches the first $n$ covariance lags is the solution of the following convex program \cite{burg1975maximum}:
\eqn \label{ME_PROBLEM}&\hat\Phi_m^\circ=\underset{\substack{\Phi_m\in\Stp}}{\mathrm{arg}\max}&  \int \log\det \Phi_m \nn\\
   &\hbox{ subject to }&   \int \Delta \Phi_m = \hat R\eeqn
The matrix $\hat R:= \left[
                 \begin{array}{cccc}
                   \hat R_0 & \hat R_1 & \ldots & \hat R_n \\
                 \end{array}
               \right]\in\Mb_{m,n}$ satisfies $\Td(\hat R)\succ 0$, \cite{STOICA_MOSES_SPECTRAL_1997}. $\hat \Phi^\circ_m$ is usually referred to as {\em maximum-entropy covariance extension}.
Because $(\hat\Phi^\circ_m)^{-1}\succ 0$ belongs to $\Qc_{m,n}$, it admits the spectral factorization
$\hat \Phi_m^\circ=\Gamma\Gamma^*$ where $\Gamma=(A\Delta^*)^{-1}$, $A\in\Rs^{m\times m(n+1)}$, is a shaping filter for the estimated process, $\hat x^\circ_m$, \cite{KAILATH_LIN_EST}. This means
that $\hat x^\circ_m$ is the output of $\Gamma$ fed by white gaussian noise (WGN), say $e$, with zero mean and variance equal to the identity:
\eq \hat x_m^\circ(t)=\sum_{j=0}^n A_j \hat x_m^\circ(t-j)+e(t)\eeq
therefore the maximum entropy estimate is an autoregressive process. In \cite{BYRNES_GUSEV_LINDQUIST_RATIONAL_COV_EXT,BYRNES_GUSEV_LINDQUIST_RATIONAL_COV_EXT_2001} it has been shown that the dual of (\ref{ME_PROBLEM})
is \eqn \label{DUAL_ME_PROBLEM}&\underset{\substack{\Phi_m^{-1}\in\Qc_{m,n}}}{\min}&  \int \left(-\log \det\Phi_m^{-1}+\Sp{\Phi_m^{-1}}{\hat \Phi_m} \right)\nn\\
   &\hbox{ subject to }&   \Phi_m\succ  0\eeqn
where $\hat \Phi_m(e^{i\vartheta}):=\sum_{j=-n}^n e^{-ij \vartheta } \hat R_j$ is the $n$-length windowed {\em correlogram} of $x^m$, \cite{STOICA_MOSES_SPECTRAL_1997}. Note that $\hat \Phi_m$ is not necessarily positive semi-definite on the unit circle.

\subsection{Spectral Density of Latent-variable Graphical Models}\label{sec_graphical_models}
We consider a real, zero-mean, stationary, purely nondeterministic, gaussian process
   $x=\{x(t);\; t\in\Zs\}$ with $m$ manifest variables and $l$ latent variables, that is $x:=\left[
                                                                                             \begin{array}{cc}
                                                                                               (x^m)^T & (x^l)^T \\
                                                                                             \end{array}
                                                                                           \right]^T$ where $x^m:=  \left[
                                                                                                                   \begin{array}{ccc}
                                                                                                                     x_1 & \ldots & x_m \\
                                                                                                                   \end{array}
                                                                                                                 \right]^T$ and
$x^l:=  \left[
                                                                                                                   \begin{array}{ccc}
                                                                                                                     x_{m+1} & \ldots & x_{m+l} \\
                                                                                                                   \end{array}
                                                                                                                 \right]^T$. Let $I\subset V_{m+l}$ be an arbitrary index set. We denote as
                                                                                                                 \eq \Xc_I=\overline{\mathrm{span}}\{x_j(t) \hbox{ s.t. } j\in I, \;t\in\Zs\}\eeq
the closure in $L_2^{m+l}(\Omega,\Ac,P)$ of the vector space of all finite linear combinations (with real coefficients) of $x_j(t)$ with $j\in I$ and $t\in\Zs$, \cite[page 3]{LINDQUIST_PICCI}.
The shorthand notation \eq \label{def_cond_independence}\Xc_{\{k\}} \botc \Xc_{\{h\}} \; | \;\Xc_{V_{m+l}\setminus \{k,h\}} \eeq means that
$\Xc_{\{k\}}$ and $\Xc_{\{h\}}$
   are conditionally independent given $\Xc_{V_{m+l}\setminus \{k,h\}}$, see \cite{ARMA_GRAPH_AVVENTI}. Therefore, (\ref{def_cond_independence}) signifies that $x_k$ and $x_h$ are conditional independent given the space linearly generated by $x_j$ with $j\in V_{m+l}\setminus \{k,h\}$. Conditional dependence relations among the variables of the process $x$ define an {\em interaction graph} $\Gc=(V_{m+l},E_{m+l})$
whose nodes represent the variables $x_1,x_2,\ldots,x_{m+l}$ and edges represent conditional dependence:
\eq (k,h)\notin E_{m+l} \;\Longleftrightarrow \; k\neq h,\; \Xc_{\{k\}} \botc \Xc_{\{h\}} \; | \;\Xc_{V_{m+l}\setminus \{k,h\}}. \eeq
The graph $\Gc$ leads to a {\em latent-variable graphical model} of the gaussian process. It admits the two layer structure illustrated in
Figure \ref{fig_ex_grah_models}: \begin{figure}[htbp]
\begin{center}
\includegraphics[width=7cm]{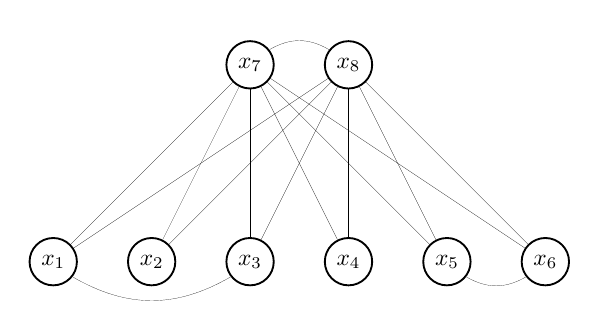}
\end{center}
\caption{Example of a latent-variable graphical model: $x_1,x_2,\ldots x_6$ are manifest variables $x_7,x_8$ are latent variables.}\label{fig_ex_grah_models}
\end{figure} latent nodes are in the upper level, and manifest nodes are in the lower level.

The graphical structure of $x$ translates into a particular decomposition of its spectral density $\Phi\in \Sph$. Starting form the block decomposition
\eq\label{Phi_p+h_inverse} \Phi =\left[
                                                  \begin{array}{cc}
                                                    \Phi_m  & \Phi_{lm} ^* \\
                                                    \Phi_{lm}  & \Phi_l  \\
                                                  \end{array}
                                                \right],\;\;
\Phi^{-1}=\left[
                                    \begin{array}{cc}
                                      \Upsilon_m & \Upsilon_{lm} ^* \\
                                      \Upsilon_{lm}  & \Upsilon_l  \\
                                    \end{array}
                                  \right]\eeq
we obtain the relationship
\eq \label{relation_complement_Schur}\Phi_m^{-1}=\Upsilon_m-\Upsilon_{lm}^*\Upsilon_l^{-1}\Upsilon_{lm}.\eeq
where we used the {\em Schur complement} pointwise.

Our main modeling assumption are that $l\leq m$ and the conditional dependence relations among the manifest variables are mostly through this limited number of latent variables. This means that the corresponding graphical model $\Gc$ has few edges between the manifest nodes, and few latent nodes. This leads to a S+L
structure for (\ref{relation_complement_Schur}), that is,
\eq \label{decmposizione_SL_assaggio}\Phi_m^{-1}=\Sigma-\Lambda, \;\; \Lambda\succeq 0 \eeq
where $\Sigma\in\Qc_{m,n}$ is sparse  and $\Lambda\in\Qc_{m,n}$ is low-rank. This means that
the support of $\Sigma$, denoted by $E_m$, contains few elements, and there exists $G\in\Rs^{l\times m(n+1)}$ with $l \ll m$ and full row rank such that $\Lambda=\Delta G^T G \Delta^*$. Accordingly, $\Phi_m^{-1}$
may be decomposed into the following two finite dimensional vector subspaces
\eqn \label{support_spaces} \Vc_{E_m}&:=& \{\Sigma\in\Qc_{m,n} \hbox{ s.t. } \Pd_{E_m^c}(\Sigma)=0\}\nn\\
\Vc_{G}&:=&\{ \Delta G^T H G \Delta^*\hbox{ s.t. } H\in\Qb_l\}.\eeqn
The sparsity of $\Sigma$ reflects the presence of few edges among the manifest nodes of $\Gc$
because of the relationship
 \eqn \label{condition_inv_spectrum} && (\Phi(e^{i\vartheta})^{-1})_{kh}=0,\;\; \forall \vartheta\in[-\pi,\pi]\;\Leftrightarrow \nn\\ &&\hspace{0.5cm}  \Xc_{\{k\}} \botc  \Xc_{\{h\}}\; |\; \Xc_{V_{m+l}\setminus\{k,h\}}\eeqn
which has been shown in \cite{ARMA_GRAPH_AVVENTI}, see also \cite{ID_DAHLHAUS,REMARKS_BRILLINGER_1996}.
The nonzero entries of $\Sigma$ therefore correspond to the (few) conditional dependence relations among the manifest variables. Accordingly, the more $\Sigma$ sparse is, the less conditional dependence relations among the manifest variables we have. Since $l\leq m$, the rank of
$\Lambda=\Upsilon_{lm}^* \Upsilon_l^{-1} \Upsilon_{lm}$ coincides with $l$, that is the number of latent variables. Accordingly the more low-rank $\Lambda$ is, the less latent variables we have. It is worth noting that (\ref{relation_complement_Schur})
is a dynamical generalization of the static decomposition \eq \label{relation_complement_Schur_static}R_m^{-1}=K_m-K_{lm}^*K_l^{-1}K_{lm}\eeq
for a zero mean gaussian random vector $x=\left[
                                      \begin{array}{cc}
                                        (x^m)^T & (x^l)^T \\
                                      \end{array}
                                    \right]\sim \Nc(0,R)$ with

\eq R=\left[
                                                  \begin{array}{cc}
                                                    R_m & R_{lm}^T \\
                                                    R_{lm} & R_l \\
                                                  \end{array}
                                                \right],\;\;
R^{-1}=\left[
                                    \begin{array}{cc}
                                      K_m & K_{lm}^T\\
                                      K_{lm} & K_l \\
                                    \end{array}
                                  \right],\eeq
see \cite{Chandrasekaran_latentvariable}. Finally, in the case $\Sigma$ is diagonal the S+L model
(\ref{decmposizione_SL_assaggio}) can be understood as a factor analysis model, \cite{deistler2007}, because conditional dependence relations among the manifest variables are only through the latent variables (or factors).

\subsection{AR Identification of Latent-variable Graphical Models}\label{sec_problem_formulation}

Let $x:=\left[
          \begin{array}{cc}
            (x^m)^T & (x^l)^T \\
          \end{array}
        \right]^T$ be an autoregressive process. We assume that a finite-length realization of $x^m$ is available, i.e. $\mathrm{x}^m(1),\mathrm{x}^m(2),\ldots \mathrm{x}^m(N)$. Regarding $x^l$, we have no data originated from it and its dimension $l$ is not even known. We would compute an estimate of the spectral density $\Phi$
of $x$. From the data, we can compute the $n$-length windowed correlogram $\hat \Phi_m$ of $x^m$. Then, the idea is to
solve the optimization problem (\ref{DUAL_ME_PROBLEM}) for the spectral density $\Phi_m$ of $x^m$ under the structural assumption (\ref{decmposizione_SL_assaggio}), but not  knowing in advance the supporting subspaces (\ref{support_spaces}). This leads us to estimate $\Vc_{E_m}$ and $\Vc_G$ first, and then estimate $\Phi_m$ consistently with the identified vector subspaces. Since the resulting estimate of the spectral density of $x^m$ obeys to (\ref{Phi_p+h_inverse}), it is then possible to recover the spectral density $\Phi$ through (\ref{Phi_p+h_inverse}) and (\ref{relation_complement_Schur}).

{\em S+L Subspace estimation}: We propose to estimate the subspaces (\ref{support_spaces}) by solving a regularized version of (\ref{DUAL_ME_PROBLEM}), that is,
 \eqn \label{PB:dual_originalpar_regularized}(\tilde \Sigma,\tilde\Lambda)=&\underset{\substack{\Sigma,\Lambda\in\Qc_{m,n}}}{\arg\min}&  \int \left(-\log (\Sigma-\Lambda)+\Sp{\Sigma-\Lambda}{\hat \Phi_m}\right)\nn\\ &&\hspace{0.2cm}+\lambda \left(\gamma\phi_1(\Sigma)+\phi_*(\Lambda)\right) \nn\\
   &\hbox{ subject to }&   \Sigma-\Lambda\succ  0\nn\\
    && \Lambda\succeq 0\eeqn
Here, $\lambda>0$ and the regularizer is a combination of two penalty functions $\phi_1$ and $\phi_*$ inducing sparsity and low-rank on $\Sigma$ and $\Lambda$, respectively. The balance between the two regularizers is tuned by $\gamma>0$. Since $(\tilde \Sigma-\tilde\Lambda)^{-1}$ represents a regularized estimate of $\Phi_m$, $\Vc_{E_m}$ is given by the support of $\tilde\Sigma$ and $\Vc_G$ by $\tilde \Lambda=\Delta G^T G \Delta^*$. Note that, for $n=0$, $\Sigma$, $\Lambda$ and $\hat\Phi_m$ are matrices, i.e. the model reduces to a gaussian random vector. In this particular situation, (\ref{PB:dual_originalpar_regularized}) boils down to the regularization problem studied in \cite{Chandrasekaran_latentvariable} for gaussian random vectors with latent variables: in that case, $\phi_1(\Sigma)$ is the $\ell_1$-norm of $\Sigma$ and $\phi_*(\Lambda)$ the nuclear norm of $\Lambda$.

{\em AR model identification}: For a fixed graphical model structure, that is, once the subspaces $\Vc_{E_m}$ and $\Vc_{G}$ have been identified, the optimal AR model is the solution to (\ref{DUAL_ME_PROBLEM}), which becomes
\eqn \label{DUAL_ME_pen}(\Sigma^\circ,\Lambda^\circ)=&\underset{\substack{\Sigma,\Lambda\in\Qc_{m,n}}}{\arg\min}&  \int \left(-\log (\Sigma-\Lambda)+\Sp{\Sigma-\Lambda}{\hat \Phi_m}\right) \nn\\
   &\hbox{ subject to }&   \Sigma-\Lambda\succ  0\nn\\
    && \Lambda \succeq 0\nn\\
&& \Sigma\in\Vc_{E_m}\nn\\
    && \Lambda\in\Vc_G  \eeqn
and the optimal estimate of $\Phi_m$ is $\hat \Phi_m^\circ=(\Sigma^\circ-\Lambda^\circ)^{-1}$.

Because the identified subspaces $\Vc_{E_m}$ and $\Vc_{G}$ depend on the regularization parameters,
a general identification procedure is as follows:
\begin{description}
\item[i)] Estimate the first $n$ covariance lags of the manifest process as in (\ref{estimated_covlag})
\item[ii)] For each $(\lambda_k,\gamma_k)$ in a given regularization path  $\{(\lambda_k,\gamma_k)\}_{k=1}^M$:
\begin{itemize}
  \item Estimate the vector subspaces $\Vc_{E_m}$ and $\Vc_G$
  \item  Compute an AR estimate $\hat \Phi_m^\circ$ of $\Phi_m$ such that $(\hat \Phi^\circ_m)^{-1}\in\Vc_{E_m}+\Vc_G$
\end{itemize}
\item[iii)] Score the identified models through a function that trades off the adherence to the data and the complexity of the models
and choose the model with the minimum score
\item[iv)] From the chosen optimal solution $\hat \Phi_m^\circ=(\Sigma^\circ-\Lambda^\circ)^{-1}$, an estimate of $\Phi$ is \eq \label{Phi_p+h_inverse_estimated}\hat \Phi= \left[
                                    \begin{array}{cc}
                                      \hat \Upsilon_m & \hat \Upsilon_{lm} ^* \\
                                      \hat\Upsilon_{lm}  & \hat\Upsilon_l  \\
                                    \end{array}
                                  \right]^{-1}\eeq where $\hat \Upsilon_m =\Sigma^\circ$, $\hat\Upsilon_{lm}$ and  $\hat\Upsilon_l$ are such that $\Lambda^\circ= \hat \Upsilon_{lm}^* \hat \Upsilon_l^{-1} \hat\Upsilon_{lm}$.  \end{description}
                                  It is worth noting that given $\Lambda^\circ$, $\hat\Upsilon_{lm}$ and $\hat\Upsilon_{l}$ are known up to an $l\times l$ invertible function. However, this is not an issue because the aim of latent variables is to explain manifest variables.
\rem Since $(\tilde \Sigma-\tilde\Lambda)^{-1}$ represents a regularized estimate of $\Phi_m$, one would wonder why it is required to solve the second problem in order to recover an estimate of $\Phi_m$. As we will see in Section \ref{sec_cov_ext}, $\hat\Phi^\circ_m$ is the maximum entropy solution of a covariance extension problem. Besides such meaningful interpretation, $\hat\Phi_m^\circ$ matches equality and inequality constraints imposed by the estimates $\hat R_j$'s which are reliable because typically we have $n\ll N$, whereas $(\tilde \Sigma-\tilde\Lambda)^{-1}$ does not. \erem

The remainder of the paper is organized as follows: the optimization problem (\ref{PB:dual_originalpar_regularized}), leading to the estimation of the sparsity
and low-rank subspaces $\Vc_{E_m}$  and $\Vc_{G}$, respectively, is studied in Section \ref{sec_top_sel}. The optimization problem (\ref{DUAL_ME_pen}), leading to the AR  model for a fixed graphical model structure, is studied in Section \ref{sec_cov_ext}. Finally, Section \ref{sec_id_procedure} provides an illustration of the full identification procedure.

\section{S+L Subspace Estimation}\label{sec_top_sel}

\subsection{Primal formulation}\label{subsec_primal_formulation}
A matrix formulation of the program (\ref{PB:dual_originalpar_regularized}) uses
(\ref{set_Qc_mn_reparametrized}), which allows to parametrize $\Sigma-\Lambda$ and $\Lambda\in\Qc_{m,n}$ as
\eqn \Sigma-\Lambda&=& \Delta X\Delta^*\in \Qc_{m,n}\nn\\
\Lambda &=&\Delta L \Delta^* \in \Qc_{m,n}\eeqn
where $X$ and $L$ are now matrix variables in the vector space $\Qb_{m(n+1)}$. Note that $\Sigma = \Delta(X+L) \Delta^*$.
Next we reformulate (\ref{PB:dual_originalpar_regularized}) in terms of $X$ and $L$.

\subsubsection{Positivity constraints $\Sigma-\Lambda\succ 0$ and $\Lambda\succeq 0$}
\lem \label{lemma_positivity}Let $\Lambda\in\Qc_{m,n}$. Then $\Lambda \succeq 0$ if and only if there exists $L\in\Qb_{m(n+1)}$ such that $L\succeq 0$.\elem
The proof is provided in Appendix \ref{proof_lemma_positivity}.

In view of Lemma \ref{lemma_positivity}, we replace the condition $\Lambda\succeq 0$ with $L\succeq 0$ and $\Sigma-\Lambda\succ 0$ with $X\succeq0$. The latter only guarantees
that $\Sigma-\Lambda\succeq 0$. However, we will show that
$X\succeq 0$ is sufficient to guarantee that $\Sigma-\Lambda\succ 0$ at the optimum of (\ref{PB:dual_originalpar_regularized}).

\subsubsection{The objective function}
Since $\Sigma-\Lambda=\Delta X\Delta^*$ with $X\succeq 0$, then there exists $A\in\Rs^{m \times m(n+1)}$
such that $X=A^TA$. By using {\em Jensen's formula}, \cite[p. 184]{AHLFORS_1953}, we obtain
\eqn \int \log\det(\Sigma-\Lambda)&=& \int \log\det(\Delta A^T A \Delta^*)\nn\\  &=& \log \det (A_0^T A_0)=\log\det X_{00}.\nn \eeqn
Clearly, the relation above holds provided that $X_{00}\succ 0$.
Moreover,
\eqn \Sp{ \Sigma-\Lambda}{\hat \Phi_m}&=&\Sp{\Delta  X \Delta^*}{\hat \Phi_m}\nn\\ &=&\Sp{\int \Delta^* \hat \Phi_m \Delta}{X}=\Sp{\Td(\hat R)}{X} \nn\eeqn
where we exploited the fact that
\eq \label{relation_R_e_correlogram}\int \Delta^* \hat \Phi_m \Delta =\Td(\hat R).\eeq
We conclude that the objective function of (\ref{PB:dual_originalpar_regularized}) admits the matrix formulation
\eqn  \label{objective_function_in_X} && \int \left(-\log\det(\Sigma-\Lambda)+\Sp{\Sigma-\Lambda}{\hat \Phi_m}\right)\nn\\ && \hspace{0.5cm}=-\log\det X_{00}+\Sp{\Td(\hat R)}{X}.\eeqn

\subsubsection{The sparsity regularizer}
Let $\Sigma\in\Qc_{m,n}$ be such that $\Sigma(e^{i\vartheta})=\sum_{j=-n}^n e^{-ij \vartheta} S_j$. Then, \eq \Pd_{E_m^c}(\Sigma)=0 \; \Longleftrightarrow \; \Pd_{E_m^c}(S_j)=0 \; j=0\ldots n.\eeq
Recall that $\Sigma=\Delta(X+L)\Delta^*$. In view of (\ref{new_reparametrization}), we obtain
\eq \label{relation_sparsity_on_D} \Pd_{E_m^c}(\Sigma)=0 \; \Longleftrightarrow \; \Pd_{E_m^c}(\Dd(X+L))=0.\eeq
We conclude that the sparsity regularizer must induce the same sparsity on the matrices $Y_j:=\Dd_j(X+L)$ with $j=0\ldots n$. In \cite{SONGSIRI_TOP_SEL_2010}, the following regularizer for $Y\in\Mb_{m,n}$ has been proposed:
\eq h_{\infty}(Y)=\sum_{k>h}\max\left\{ |(Y_0)_{hk}|,\max_{j=1\ldots n}|(Y_j)_{hk}|,\max_{j=1\ldots n} |(Y_j)_{kh}|\right\}.\eeq
Let $v_{kh}$, with $k>h$, be the vector of $(k,h)$ and $(h,k)$ entries of the coefficients $Y_j$ with $j=0\ldots n$. Therefore, \eq h_\infty(Y)=\sum_{k>h} \| v_{kh}\|_\infty\eeq
where $\|\cdot \|_\infty$ denotes the $\ell_\infty$-norm. On the other hand, $h_\infty(Y)$ is the $\ell_1$-norm of the vector having (nonnegative) entries $\| v_{kh}\|_\infty$ with $k>h$. Accordingly, $h_\infty(Y)$ encourages sparsity among $v_{kh}$'s, that is induces the same sparsity on the matrices $Y_j$ $j=0\ldots n$.
\subsubsection{The low-rank regularizer}
\prop \label{prop_cov_envelope} Given $\Lambda\in\Ac_m$, we define the convex function
\eq \label{def_low_rank_regularizer}\phi_*(\Lambda):=\sum_{k=1}^{m} \int \sigma_k(\Lambda)\eeq
and the restricted rank function
\eq \mathrm{rank}^\prime(\Lambda):=\left\{
  \begin{array}{ll}
    \mathrm{rank}(\Lambda), & \|\Lambda\|\leq 1 \\
    +\infty, & \hbox{otherwise.}
  \end{array}
\right.\eeq
 Then, the convex hull of $\mathrm{rank}^\prime(\Lambda)$ is \eq \label{def_convex_hull} \left\{
  \begin{array}{ll}
    \phi_*(\Lambda), & \|\Lambda\|\leq 1 \\
    +\infty, & \hbox{otherwise.}
  \end{array}
\right.\eeq \eprop
The proof is provided in Appendix \ref{proof_conv_envelope}.

We conclude that $\phi_*(\Lambda)$ defined in (\ref{def_low_rank_regularizer}) is the adequate regularizer of $\mathrm{rank}(\Lambda)$.
Since $\Lambda\succeq 0$, $\sigma_k(\Lambda(e^{i\vartheta}))$ represents the $k$-th eigenvalue of $\Lambda(e^{i\vartheta})$. Thus, $\phi_*(\Lambda)=\tr\int \Lambda$.
Finally, \eq    \phi_*(\Lambda)=\tr\int \Delta L \Delta^* =\tr\left( L\int \Delta^*\Delta\right)=\tr(L)\nn\eeq
where we exploited the fact that \eq \label{proprieta_int_exp}\int e^{i j\vartheta }=\left\{
                                                                                       \begin{array}{ll}
                                                                                         1,& j=0 \\
                                                                                         0,& j\neq 0.
                                                                                       \end{array}
                                                                                     \right.
\eeq

\subsubsection{Primal Formulation}
By collecting the results in 1)-4), we rewrite (\ref{PB:dual_originalpar_regularized}) as
\eqn \label{PB:dual_Delta_regularized}(X^\circ,L^\circ)=&\hspace{-0.5cm}\underset{\substack{X ,L\in\Qb_{m(n+1)}}}{\arg\min}& \hspace{-0.5cm}- \log\det X_{00}+\Sp{\Td(\hat R)}{X}\nn\\ && \hspace{-0.1cm}+\lambda\gamma h_\infty(\Dd(X+L)) +\lambda\tr(L)\nn\\
   & \hspace{-0.5cm}\hbox{ subject to }&    \hspace{-0.5cm}X_{00}\succ 0,\;X\succeq 0,\; L\succeq 0\eeqn
 Formulation (\ref{PB:dual_Delta_regularized}) and (\ref{PB:dual_originalpar_regularized}) are equivalent provided that
$\Delta X^\circ \Delta^* \succ 0$. Finally, it is worth noting that (\ref{PB:dual_Delta_regularized}) is a generalization of the regularized problem studied in \cite{SONGSIRI_TOP_SEL_2010}. The problem formulations coincide when $L=0$, that is for estimating an AR process having a sparse graphical model but no latent variables.

\subsection{Dual formulation}

 We show that (\ref{PB:dual_Delta_regularized}) does admit a solution by exploiting duality theory. First, note that (\ref{PB:dual_Delta_regularized}) is strictly feasible (pick $X=I$ and $L=I$), thus {\em Slater's condition} holds. Accordingly, the duality gap between (\ref{PB:dual_Delta_regularized}) and its dual problem is equal to zero.
We introduce a new variable $Y\in\Mb_{m,n}$ in (\ref{PB:dual_Delta_regularized}) to obtain the following equivalent problem
\eqn \underset{\tiny \substack{X ,L\in\Qb_{m(n+1)}\\ Y\in\Mb_{m,n}}}{\arg\min}& \hspace{-0.5cm}- \log\det X_{00}+\Sp{\Td(\hat R)}{X}+\lambda\gamma h_\infty(Y) +\lambda\tr(L)\nn\\
   \hbox{ subject to }& \hspace{-3.5cm}  X_{00}\succ 0,\;X\succeq 0,\; L\succeq 0\nn\\
   & \hspace{-4.8cm} Y=\Dd(X+L)\nn\eeqn
The Lagrangian is \eqn  && \hspace{-0.6cm}\Lc(X,L,Y,U,V,Z)\nn\\
&&\hspace{-0.4cm} =-\log\det X_{00}+\Sp{\Td(\hat R)}{X}+\lambda\gamma h_\infty(Y)+\lambda \tr(L)\nn\\
&& -\Sp{U}{X}-\Sp{V}{ L}+ \Sp{Z}{\Dd(X+L)-Y}\nn \\
&&\hspace{-0.4cm}=-\log\det X_{00}+\Sp{\Td(\hat R)-U}{X}+\Sp{\lambda I -V}{ L}\nn\\
&&+\lambda\gamma h_\infty(Y)-\Sp{Z}{Y}+\Sp{\Td(Z)}{X+L} \nn \\
&&\hspace{-0.4cm} =-\log\det X_{00}+\Sp{\Td(\hat R)+\Td(Z)-U}{X}\nn\\
&&+ \Sp{\lambda I +\Td(Z)-V}{ L}+\lambda\gamma h_\infty(Y) -\Sp{Z}{Y} \nn  \eeqn
where $U,V\in\Qb_{m(n+1)}$ such that $U,V\succeq 0$ and $Z\in\Mb_{m,n}$.
The dual function is the infimum of $\Lc$ over $X$, $L$ and $Y$.  We start by minimizing with respect to $Y$. The Lagrangian $\Lc$ depends on $Y$
only through the term
\eq \label{term_in_Sb}\lambda\gamma h_\infty(Y)-\Sp{Z}{Y}, \eeq
which was shown, \cite{SONGSIRI_TOP_SEL_2010}, to be bounded below only if \eqn \label{cond_S1} && \mathrm{diag}(Z_j)=0,\;\;j=0\ldots n\\ \label{cond_S2} && \sum_{j=0}^n |(Z_j)_{kh}|+|(Z_j)_{hk}|\leq \lambda\gamma,\;\; k\neq h,\eeqn
in which case the infimum is equal to zero. The partial minimization of the Lagrangian over $Y$ is therefore \eqn && \hspace{-0.6cm}\inf_{Y} \Lc =\left\{  \begin{array}{ll}
 -\log\det X_{00}+\Sp{\Td(\hat R)+\Td(Z)-U}{X} &\\
\hspace{0.5cm} +\Sp{\lambda I +\Td(Z)-V}{ L} & \hspace{-0.6cm}\hbox{(\ref{cond_S1}),\; (\ref{cond_S2})} \\
                                                                -\infty & \hspace{-0.6cm}\hbox{otherwise.}
                                                              \end{array}
                                                        \right.\nn
\eeqn
Likewise, the Lagrangian $\Lc$ depends on $L$ only through the term $\Sp{\lambda I+\mathrm{T}(Z)-V}{L}$, which is
bounded below only if
\eq \label{cond_L}\lambda I+\mathrm{T}(Z)-V=0, \eeq in which case the infimum is equal to zero.
Thus,
\eqn  &&\hspace{-0.7cm}\inf_{L,Y}  \Lc =\left\{
                                                              \begin{array}{ll}
 \hspace{-0.2cm}-\log\det X_{00}+\Sp{\Td(\hat R)+\Td(Z)-U}{X} & \hspace{-0.2cm}\hbox{(\ref{cond_S1})-(\ref{cond_L})} \\
      \hspace{-0.2cm}                                                          -\infty & \hspace{-0.4cm}\hbox{otherwise.}
                                                              \end{array}
                                                            \right. \nn
\eeqn Finally, the terms in $X_{00}$ are bounded below if and only if \eq \label{cond_delta2}(\mathrm{T}(Z)+\Td(\hat R)-U)_{00}\succ 0\eeq and if (\ref{cond_delta2}) holds, they are minimized by $X_{00}=(\mathrm{T}(Z)+\Td(\hat R)-U)_{00}^{-1}$. The Lagrangian is linear in the remaining variables $X_{kh}$, and therefore bounded below (and identically zero) only if \eq \label{cond_delta} (\mathrm{T}(Z)+\Td(\hat R)-U)_{kh}=0 \;\;\forall \;(k,h)\neq (0,0).\eeq The final expression for the dual functional is
\eq \label{funzionale_duale_regolarizzato}\inf_{X,L,Y} \Lc=\left\{
                                                              \begin{array}{ll}
                                                                \log\det (\mathrm{T}(Z)+\Td(\hat R)-U)_{00}+m & \hbox{(\ref{cond_S1})-(\ref{cond_delta})} \\
                                                                -\infty & \hbox{otherwise.}
                                                              \end{array}
                                                            \right. \eeq The dual problem consists in maximizing the dual functional (\ref{funzionale_duale_regolarizzato})  with respect to $U$, $V$ and $Z$ subject to the constraints $U\succeq 0$ and $V\succeq 0$. Moreover, eliminating the slack variables $U$ and $V$, and adding the variable $W:=(\mathrm{T}(Z)+\Td(\hat R)-U)_{00}$ the dual problem takes the final form
\eqn \label{PB:dual_dual_Delta_regularized}&\underset{\substack{W\in\Qb_m\\ Z\in\Mb_{m,n}}}{\max} & \log\det W +m\nn\\ &\hbox{subject to }&  W\succ 0\nn\\
&&  \Td(\hat R)+\mathrm{T}(Z)\succeq \left[
                                 \begin{array}{cc}
                                   W & 0 \\
                                   0 & 0 \\
                                 \end{array}
                               \right]\nn\\
&& \mathrm{diag}(Z_j)=0,\;\;j=0\ldots n\nn \\ &&
\sum_{j=0}^n |(Z_j)_{kh}|+|(Z_j)_{hk}|\leq \lambda\gamma,\;\; k\neq h\nn\\
&& \lambda I+\mathrm{T}(Z)\succeq 0   \eeqn
 \prop \label{prop_dual_dual_Delta_regularized}Problem (\ref{PB:dual_dual_Delta_regularized}) admits a solution. \eprop
The proof is provided in Appendix \ref{proof_prop_dual_dual_Delta_regularized}.

From the next statement we conclude that Problem (\ref{PB:dual_originalpar_regularized}) admits a solution.
\prop \label{prop_equivalence}Problem (\ref{PB:dual_Delta_regularized}) admits a solution $(X^\circ,L^\circ)$  such that $\Delta X^\circ \Delta^*\succ 0$. Accordingly   (\ref{PB:dual_originalpar_regularized}) and (\ref{PB:dual_Delta_regularized}) are equivalent. Moreover, $X^\circ $ is unique.\eprop
The proof is provided in Appendix \ref{proof_prop_equivalence}.

It is worth noting that (\ref{PB:dual_dual_Delta_regularized}) is easier to solve than (\ref{PB:dual_Delta_regularized}), because the objective function in (\ref{PB:dual_dual_Delta_regularized}) is smooth.

\subsection{Estimation of the Vector Subspaces}

The vector subspace $\Vc_{E_m}$ is given by the support of $\tilde \Sigma=\Delta (X^\circ+L^\circ) \Delta^*$. In view of (\ref{new_reparametrization}), we obtain
\eq \label{recover_E_c}E_m^c=\{(k,h)\in V_m\times V_m \hbox{ s.t. } (\Dd(X^\circ+L^\circ))_{kh}=0\}\eeq
and hence also $\Vc_{E_m}$.  Since $\tilde \Lambda=\Delta L^\circ \Delta^*$, the vector subspace $\Vc_G$ is the column space of $L^\circ$, given by the decomposition $L^\circ =G^T G$ where $G$ is a full row rank matrix.

Next, we show how to recover $(X^\circ,L^\circ)$ from an optimal solution $(W^\circ,Z^\circ)$
of the smooth convex optimization program (\ref{PB:dual_dual_Delta_regularized}). Such a recovering scheme also provides sufficient conditions
for the uniqueness of the two vector subspaces. Regarding $X^\circ$, let $B\in\Rs^{m\times m(n+1)}$ the solution of the {\em Yule-Walker}
equation \eq  (\Td(\hat R)+\Td(Z^\circ))B^T=\left[
                                              \begin{array}{c}
                                                W^\circ \\
                                                0 \\
                                              \end{array}
                                            \right],\;\; B_0=I
\eeq then $X^\circ=B^T (W^\circ)^{-1}B$, see Appendix \ref{proof_prop_equivalence} for more details. Next, we deal with the recovering of $L^\circ$.
 Because of the strong duality between (\ref{PB:dual_Delta_regularized}) and (\ref{PB:dual_dual_Delta_regularized}), we have
  \eq \label{dual_gap_zero_UL} \Sp{V^\circ}{L^\circ}=0\eeq where $V^\circ:=\lambda I+\mathrm{T}(Z^\circ)$, see (\ref{cond_L}). If $V^\circ$ is a full rank matrix then, in view of (\ref{dual_gap_zero_UL}), $L^\circ=0$ is the unique solution, $\Vc_G=\{0\}$ and $\Vc_{E_m}$ is univocally characterized by (\ref{recover_E_c}). Otherwise, let $l>0$ be the dimension of the nullspace of $V^\circ$. Then there exists a full row rank matrix $G\in\Rs^{l\times m(n+1)}$ such that $V^\circ G^T=0$. Since $V^\circ,L^\circ \succeq 0$, from (\ref{dual_gap_zero_UL}) it follows that \eq L^\circ=G^T H G\eeq where $H$, unknown, belongs to $\Qb_l$  and $H\succeq 0$.
Therefore, $L^\circ$ is known up to the (scaling) factor $H$. The minimization of (\ref{term_in_Sb}) under constraints (\ref{cond_S1}) and (\ref{cond_S2})
is equivalent to the minimization of the non-negative function \eqn \label{def_f_in_k_h} &&\max\{|(Y_0)_{kh}|,\max_{j=1\ldots n}|(Y_j)_{kh}|,\max_{j=1\ldots n}|(Y_j)_{hk}|\}\nn\\ && \hspace{0.5cm}\times\left(\lambda\gamma -\sum_{j=0}^n |(Z_j)_{kh}|+|(Z_j)_{hk}|\right), \eeqn
for each $k>h$, subject to the constraint that their sum is bounded by $\lambda\gamma h_\infty(Y)-\Sp{Z}{Y}$.
Since the optimal value of (\ref{term_in_Sb}) is always equal to zero, then the optimal value of (\ref{def_f_in_k_h}) is equal to zero for each $k>h$.
Thus, if $\sum_{j=0}^n |(Z_j)_{kh}|+|(Z_j)_{hk}|<\lambda \gamma$ then \eq \max\{|(Y_0)_{kh}|,\max_{j=1\ldots n}|(Y_j)_{kh}|,\max_{j=1\ldots n}|(Y_j)_{hk}|\}=0\eeq and $(Y_j)_{kh}=(Y_j)_{hk}=0$ with $j=0\ldots n$.
Since $Y=\Dd(X+L)$, $\sum_{j=0}^n |(Z_j)_{kh}|+|(Z_j)_{hk}|<\lambda \gamma$ implies that $(\Dd_j(X+L))_{kh}=(\Dd_j(X+L))_{hk}=0$ with $j=0\ldots n$.
Accordingly, $H$ is obtained by solving the following system of linear equations
\eq  \label{constr_Sigma}  (\Dd_j(X^\circ+G^T H G))_{kh}=0  \;\;j=0\ldots n , \; \forall\;(k,h)\in I.
 \eeq
where
\eq  I:=\left\{(k,h) \hbox{ s.t. } k\neq h,\;\sum_{j=0}^n |(Z_j)_{kh}|+|(Z_j)_{hk}|<\lambda\gamma \right\}.\eeq
Note that (\ref{constr_Sigma}) is a system of $(n+1)\times |I|$ equations with $l(l+1)/2$ unknowns (i.e. the number of independent parameters in $H$). For $\lambda\gamma$
and $\gamma$ sufficiently large, $|I|$ would be sufficiently large and $l$ sufficiently small, respectively, so that (\ref{constr_Sigma}) admits a unique solution.
We stress the fact that it may happen that (\ref{constr_Sigma}) has not unique solution even $l(l+1)/2\ll (n+1)\times |I| $. As observed in \cite{Chandrasekaran_latentvariable}, this is more likely when $\Vc_{G}$ contains sparse elements, that is the latent variables are not sufficiently ``diffuse''
across the manifest variables, or $\Vc_{E_m}$ contains elements with a low degree of sparsity, that is the are manifest variables conditionally dependent to too many other manifest variables. Both cases may lead to a non-identifiability of the AR model solution to Problem (\ref{DUAL_ME_pen}) because some sparse and low-rank components are not distinguishable. One avoids those situations checking that (\ref{constr_Sigma})
has unique solution. We formalize the above observation.
\prop \label{prop_transversality}If (\ref{constr_Sigma}) admits a unique solution, then $\Vc_{E_m}$ and $\Vc_G$  are unique and have {\em transverse intersection}, i.e. $\Vc_{E_m} \cap \Vc_G=\{0\} $.
\eprop
The proof is provided in Appendix \ref{proof_prop_transversality}.

The transversality condition means that any element of $\Vc_{E_m}+\Vc_G$ admits a unique decomposition into the two subspaces. We will see in Section \ref{sec_cov_ext} that this condition guarantees the uniqueness of the solution to Problem (\ref{DUAL_ME_pen}).

\section{AR Model Identification}\label{sec_cov_ext}
The convex formulation of the convex optimization Problem (\ref{DUAL_ME_pen})
parallels the developments in the previous section. We adopt the parametrization
\eqn \Sigma- \Lambda &=&\Delta X \Delta^*\nn\\
\Lambda &=&\Delta L\Delta^*=\Delta G^T H G \Delta^* \eeqn
where the matrix unknowns are $X\in\Qb_{m(n+1)}$ and $H\in\Qb_l$. Note that $\Sigma=\Delta(X+G^T H G)\Delta^*$.
The positivity conditions $\Sigma-\Lambda\succ 0$ and $\Lambda\succeq 0$ are replaced by $X\succeq 0$ and $H\succeq 0$, respectively. Also in this case $X\succeq 0$ only guarantees that $\Sigma-\Lambda\succeq 0$.
In view of (\ref{relation_sparsity_on_D}), condition $\Sigma\in\Vc_{E_m}$ is replaced by $\Pd_{E_m^c}(\Dd(X+G^THG))=0$. Clearly condition $\Lambda\in\Vc_G$ follows from the chosen parametrization.
Finally, the objective function is given by (\ref{objective_function_in_X}) provided that $X_{00}\succ 0$. The convex program (\ref{DUAL_ME_pen}) thus admits the matrix formulation
\eqn \label{PB:dual_Delta}&\underset{\substack{X\in\Qb_{m(n+1)}\\ H\in\Qb_l}}{\min} &  - \log\det X_{00}+\Sp{\Td(\hat R)}{X} \nn\\
   &\hbox{ subject to} & X_{00}\succ 0,\; X\succeq 0,\; H\succeq 0\nn\\
   &&\Pd_{E_m^c}(\Dd(X+ G^T H G))=0
\eeqn
Both formulations are equivalent provided that the optimal solution, say $(X^\circ,H^\circ)$, is such that $\Delta X^\circ\Delta^*\succ 0$.
\prop  \label{prop_existence_sol} Problem (\ref{PB:dual_Delta}) does admit a solution. Moreover, $\Delta X^\circ \Delta$ is unique and such that $\Delta X^\circ \Delta\succ 0$.\eprop
The proof is provided in Appendix \ref{proof_existence_sol}.

The optimal spectral density $\hat \Phi_m^\circ$ thus admits the matrix decomposition
\eq (\hat \Phi_m^\circ)^{-1}=\Delta X^\circ \Delta^*= \underbrace{\Delta (X^\circ+G^T H^\circ G)\Delta^*}_{\in\Vc_{E_m}}-\underbrace{\Delta G^T H^\circ G\Delta^*}_{\in\Vc_G}\eeq
which is unique when $\Vc_{E_m}$ and $\Vc_G$ have transverse intersection.

\cor The AR latent-variable graphical model solution to Problem (\ref{DUAL_ME_pen}) is unique when $\Vc_{E_m}$ and $\Vc_{G}$ are estimated from (\ref{PB:dual_originalpar_regularized}) with $\lambda\gamma$ and $\lambda$ sufficiently large.\ecor

We now give an important interpretation of the optimal solution of (\ref{DUAL_ME_pen}). Consider the following covariance extension problem.

\pb \label{prob_cov_ext}Find $\Phi_m\in\Sc_m$ such that \eqn && P_{E_m}\left(\int \Delta \Phi_m-\hat R \right)=0 \nn\\ && \int  G \Delta^*\Phi_m \Delta G^T\succeq   G \Td(\hat R) G^T.\eeqn \epb

The condition
\eq  \label{condition_first_n_covlags}\int \Delta \Phi_m=\hat R\eeq
implies that $\Pd_{E_m}\left(\int \Delta \Phi_m-\hat R\right)=0$. Moreover, (\ref{condition_first_n_covlags}) is equivalent to
$\int \Delta^* \Phi_m \Delta=\Td(\hat R)$ which implies that $\int G\Delta^* \Phi_m \Delta G^T \succeq G\Td(\hat R) G^T$. Accordingly, Problem \ref{prob_cov_ext} is a relaxation of the classic covariance extension problem.
The next theorem shows that $\hat\Phi_m^\circ$ is the maximum entropy solution of Problem \ref{prob_cov_ext}.
\teo Problem (\ref{DUAL_ME_pen}) is the dual of the convex optimization problem
\eqn \label{PB:ME_latent_graph}& \underset{\Phi_m\in\Stp}{\max} & \int \log\det \Phi_m\nn\\
   &\hbox{subject to}  & P_{E_m}\left(\int \Delta \Phi_m-\hat R \right)=0 \nn\\
   && \int  G \Delta^*\Phi_m \Delta G^T\succeq   G \Td(\hat R) G^T\eeqn
\eteo
\proof
Note that, (\ref{PB:ME_latent_graph}) is a relaxation of (\ref{ME_PROBLEM}).
Moreover, (\ref{ME_PROBLEM}) admits solution (and thus it is feasible), because $\Td(\hat R)\succ 0$. Accordingly, (\ref{PB:ME_latent_graph}) is feasible.
Moreover, we only have linear inequality constraints in (\ref{PB:ME_latent_graph}) which implies the {\em refined Slater's condition} \cite{BOYD_CONVEX_OPTIMIZATION}.  Thus we have strong duality for (\ref{PB:ME_latent_graph}) and its dual. The Lagrange
functional is: \eqn && \hspace{-0.8cm}\Lc(\Phi_m,S,H)=\int\log\det \Phi_m- \Sp{\Pd_{E_m}\left(\int \Delta \Phi_m-\hat R \right)}{S}\nn\\
&& \hspace{1.7cm}+\Sp{\int G \Delta^*\Phi_m \Delta G^T- G \Td(\hat R) G^T}{H}\eeqn
where $H\in\Qb_l$ such that $H\succeq 0$, and $S\in\Mb_{m,n}$.
Moreover,
\eqn && \hspace{-0.6cm}\Lc(\Phi_m,S,H)\nn\\ &&\hspace{-0.5cm} =\int\log\det \Phi_m- \Sp{\int \Delta\Phi_m-\hat R}{\Pd_{E_m}\left(S\right)}\nn\\
 && \hspace{0.0cm}+\Sp{\int G  \Delta^*(\Phi_m-  \hat \Phi_m) \Delta  G^T}{H}\nn\\
&&\hspace{-0.5cm}=\int\log\det \Phi_m- \Sp{\int \Delta\Phi_m-\hat R}{\Pd_{E_m}\left(S\right)}\nn\\
&& \hspace{0.0cm}+\Sp{\Phi_m -\hat \Phi_m}{\Delta G^T H G \Delta^*}\eeqn
where we exploited (\ref{relation_R_e_correlogram}) and the fact that $\Pd_{E_m}$ is a self-adjoint operator.
By defining $\Sigma:=\Pd_{E_m}(S_0)+\frac{1}{2}\sum_{j=1}^n e^{-ij\vartheta} \Pd_{E_m} (S_j)+e^{ij\vartheta} \Pd_{E_m} (S_j)^T\in\Vc_{E_m}$ and $\Lambda:=\Delta G^T H G \Delta^*\in \Vc_G$ such that $\Lambda \succeq 0$, we obtain the following compact notation for the Lagrangian
\eqn &&\hspace{-0.8cm}\Lc (\Phi_m,\Sigma,\Lambda)\nn\\&& \hspace{-0.8cm}= \int\left(\log\det\Phi_m-\Sp{\Phi_m-\hat \Phi_m}{\Sigma}+\Sp{\Phi_m-\hat \Phi_m}{\Lambda}\right).\nn\eeqn
Since $\Lc(\cdot,\Sigma,\Lambda)$ is strictly concave over $\Stp$, its unique maximum point is given by annihilating its first variation in each direction $\delta \Phi_m\in \Ld_2^{m\times m}$:
\eq \delta \Lc(\Phi_m,\Sigma,\Lambda;\delta\Phi_m)=\tr\int \left((\Phi_m^{-1}-\Sigma+\Lambda)\delta\Phi_m\right)\eeq
Note that $\Phi_m^{-1}-\Sigma+ \Lambda \in \Ld_2^{m\times m}$, thus the first variation is zero for each $\delta \Phi_m$ if and only if $\Phi_m^{-1}-\Sigma + \Lambda =0$. Accordingly, if $\Sigma-\Lambda \succ 0$ then the unique maximum point of $\Lc(\cdot,\Sigma,\Lambda)$ is \eq \label{optimal_form}\hat \Phi_m^\circ:=\left(\Sigma-\Lambda\right)^{-1}\eeq
with $\Sigma\in\Vc_{E_m}$ and $\Lambda\in\Vc_G$ such that $\Lambda\succeq 0$.
Then, by substituting (\ref{optimal_form}) in the Lagrangian we obtain, up to a constant term, the objective function of (\ref{DUAL_ME_pen}). \qed\\

The interpretation of the convex program (\ref{DUAL_ME_pen}) as the dual of a covariance extension problem is insightful. First, it coincides with the problem considered in \cite{ARMA_GRAPH_AVVENTI} for the AR case,
since the solution satisfies $G=0$ when the inequality constraint in (\ref{PB:ME_latent_graph}) is removed. On the other hand, it is worth noting that \cite{ARMA_GRAPH_AVVENTI} considers ARMA models which are more general than the AR ones. Second, both constraints in (\ref{PB:ME_latent_graph}) have a clear interpretation: the equality constraint  imposes that the optimal spectral density $\hat \Phi_m^\circ$ matches the estimated covariance lags $\hat R_0\ldots \hat R_n$ in the positions specified by $E_m$.
Regarding the inequality constraint, consider the stochastic process \eq y(t)=\sum_{j=0}^n G_j x^m(t-j)\eeq
whose variables are linear combinations of the $m$ manifest variables in a time window of length $n$. Accordingly, $y$ encodes information about $x^m$.
It is readily checked that  \eq \Es[y(t)y(t)^T]=\sum_{k,h=0}^n G_k R_{h-k} G_h^T= G \Td(R)G^T.\eeq
The inequality constraint therefore imposes that the covariance matrix of $y$ is lower bounded by the one estimated from the data, i.e. $G \Td(\hat R)G^T$.

\section{Numerical Examples}\label{sec_id_procedure}
\subsection{Synthetic Example}\label{susec_synthetic_ex}
We consider an AR latent-variable graphical model of order $n=1$ with $m=15$ manifest variables, $l=1$ latent variable. Its interaction graph is depicted in Figure
\ref{Figura1}(a). We generate a data sequence of length $N=500$ for the manifest process and we apply
the identification procedure outlined at the end of Section \ref{sec_problem_formulation}.
In Figure \ref{Figura1}(b) we depict the latent-variable graphical models obtained for different values of $\lambda$ and $\lambda\gamma$.
\begin{figure*}[htbp]
\begin{center}
\includegraphics[width=15cm]{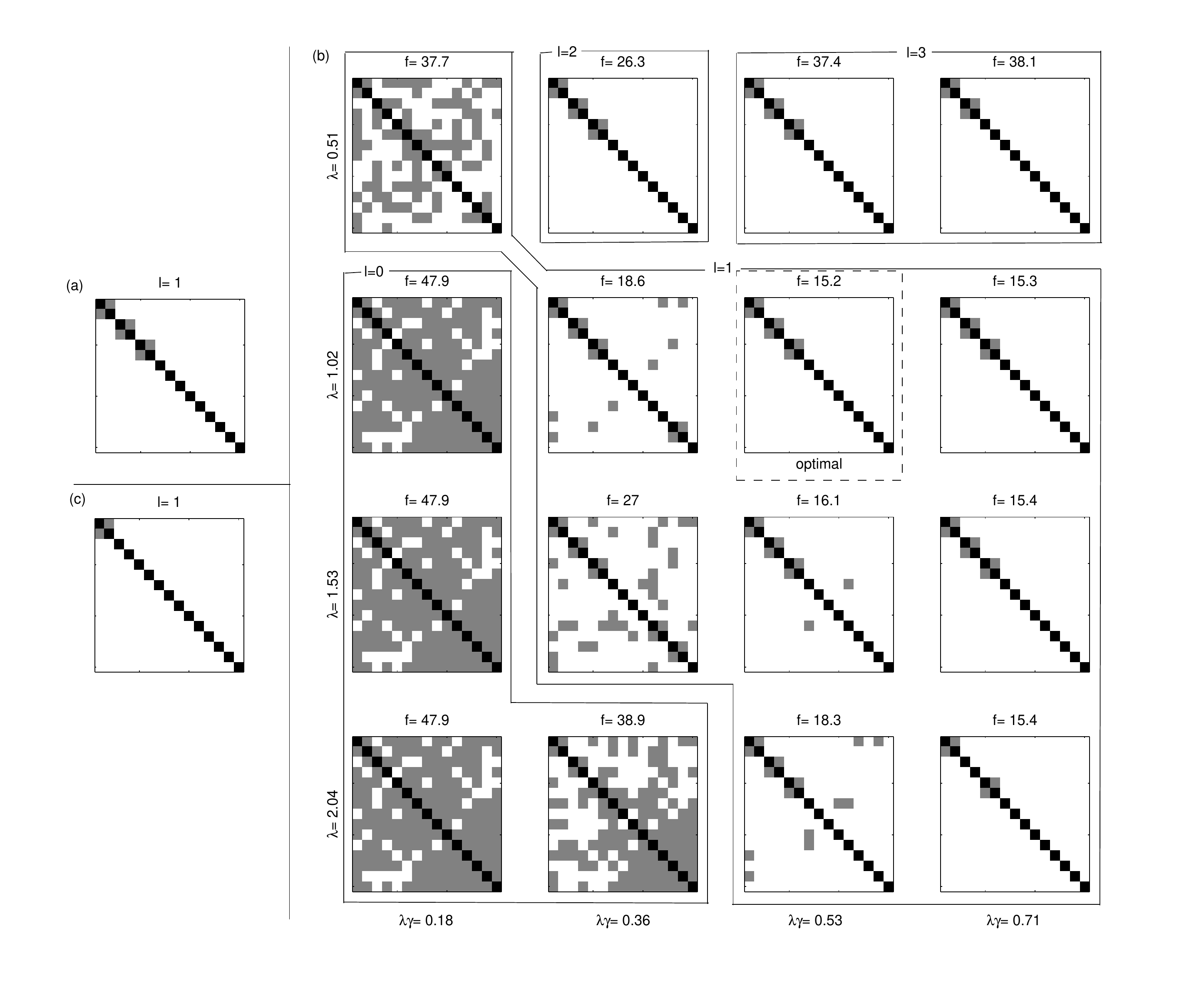}
\end{center}
 \caption{(a) Interaction graph of the generated model. (b) Interaction graphs of optimal models estimated for  $n=1$ and for different values of $\lambda$ and $\lambda\gamma$. (c)
Interaction graph of the optimal model estimated with $n=0$. Each figure shows the interaction graph for the manifest variables: grey denotes an edge, white denotes no edge, and black denotes a manifest node. The number of latent variables and the value of the score function is indicated on the top of each figure.}\label{Figura1}
\end{figure*}
\begin{figure}[htbp]
\begin{center}
\includegraphics[width=9.5cm]{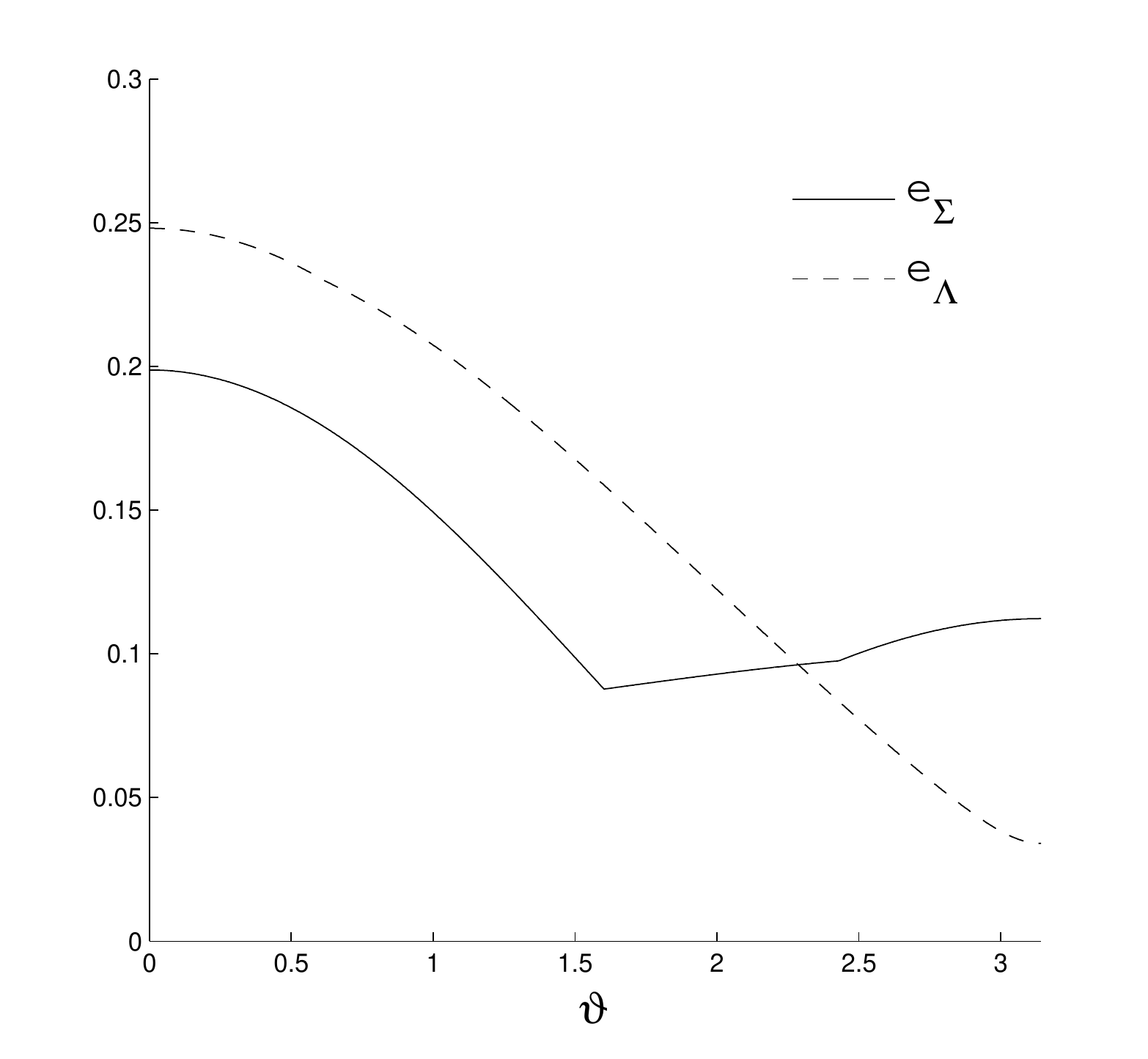}
\end{center}
 \caption{Normalized estimation errors $\mathrm{e}_\Sigma(e^{i\vartheta})=\frac{\| \Sigma(e^{i\vartheta})-\Sigma^\circ(e^{i\vartheta})\|_2}{\|\Sigma\|}$ and  $\mathrm{e}_\Lambda(e^{i\vartheta})=\frac{\| \Lambda(e^{i\vartheta})-\Lambda^\circ(e^{i\vartheta})\|_2}{\|\Lambda\|}$ as a function of $\vartheta\in[0,\pi]$ for the data set in Section \ref{sec_id_procedure}.}\label{Figura2}
\end{figure}
Not surprisingly, increasing the rank regularization parameter $\lambda$ favors few latent variables, whereas by increasing the sparsity regularization parameter $\lambda\gamma$ favors few conditional dependence relations among the manifest variables.

To discriminate among models, we consider the following score function:
\eq \label{fitness_function}f(E_{m},l,\hat \Phi_{m}^\circ,\hat \Phi_{C})=\Ds(\hat \Phi_{C}\|\hat \Phi_{m}^\circ)\times p.\eeq
Here, $\hat \Phi_C$ is the smoothed {\em correlogram} of $x^m$ computed from the data by using the {\em Bartlett window}, \cite{STOICA_MOSES_SPECTRAL_1997}.
The cost \eqn \Ds(\hat \Phi_C\| \hat \Phi_{m}^\circ):=\frac{1}{2}\left(\int
\left(\log\det(\hat \Phi_C^{-1}\hat \Phi_{m}^\circ)\right.\right.\nn\\ \left.\left.+\Sp{\hat \Phi_C}{(\hat \Phi_{m}^\circ)^{-1}}\right)-m\right)\eeqn
is the  {\em relative entropy rate}, \cite{COVER_THOMAS}, between
 $\hat \Phi_C$ and $\hat \Phi_{m}^\circ$. Thus, it ranks the adherence of $\hat \Phi_{m}^\circ$ to the data.  The term \eq p=(|E_{m}|-m)+m l\eeq
is the total number of edges in the latent-variable graphical model. Thus, $p$ places a penalty on models with high complexity.
An alternative choice for the score function would be $\Ds(\hat \Phi_{C}\|\hat \Phi_{m}^\circ)+\alpha(N) p$ where the weighting $\alpha(N)$ is the trade-off parameter
between the adherence to the data and the complexity of the model. Typically $\alpha(N)$ is a decreasing function in $N$ because the data should reveal the simple structure as $N$ increases. The authors of  \cite{SONGSIRI_GRAPH_MODEL_2010} recommend the choices $\alpha(N)=N^{-1}$ and $\alpha(N)=N/ \log N$. In contrast, the authors of \cite{ARMA_GRAPH_AVVENTI} recommend the score function (\ref{fitness_function}) because  it is robust to scaling. Based on (\ref{fitness_function}), the minimum value of $f$ is equal to 15,2 reached
with $\lambda=1.02$ and $\lambda\gamma=0.53$. Its interaction graph coincides with the true one. Figure \ref{Figura2} provides a graph of the normalized estimation errors of $\Sigma$ and $\Lambda$ at each frequency:
\eqn \mathrm{e}_\Sigma(e^{i\vartheta})&=&\frac{\| \Sigma(e^{i\vartheta})-\Sigma^\circ(e^{i\vartheta})\|_2}{\|\Sigma\|}\nn\\ \mathrm{e}_\Lambda(e^{i\vartheta})&=&\frac{\| \Lambda(e^{i\vartheta})-\Lambda^\circ(e^{i\vartheta})\|_2}{\|\Lambda\|}.\eeqn
We found similar results by varying the sample data. Finally, we applied the same identification procedure with $n=0$, i.e. by estimating a gaussian random vector. The estimated interaction graph in Figure \ref{Figura1}(c)
does not recover the generated model. This suggests the potential benefit of AR modeling in the estimation of latent-variable graphical models.
\subsection{International Stock Markets}
The data used in this simulation consists of a time series of daily stock markets indices at closing time, in terms of local currency units, of twenty-two financial markets. The twenty-two countries an their respective price indices are: Australia (All Ordinaries index denoted AU), New Zealand (50 Gross index denoted NZ), Singapore (STI index denoted SG), Hong Kong (Hang Seng index denoted HK), China (SSE Composite index denoted CH), Japan (Nikkei225 index denoted JA), Korea (KOSPI Composite index denoted KO), Taiwan (Weighted index denoted TA), Brazil (IBOVESPA index denoted BR), Mexico (IPC index denoted ME), Argentina (Merval index denoted AR), Swiss (SMI index denoted SW), Greece (Athen Composite index denoted GR), Belgium (BFX index denoted BE), Austria (ATX index denoted AS), Germany (DAX index denoted GE), France (CAC 40 index denoted FR), Netherlands (AEX index denoted NL), United Kingdom (FTSE 100 index denoted UK), Unites States (S\&P500 denoted US), Canada (S\&PTSX Composite index denoted CA) and Malaysia (KLCI index denoted MA). The data are obtained from the website at \url{http://finance.yahoo.com/}. The sample period is from $4$th January 2012 up to $31$th December 2013. For each index, we compute the return between the trading day
$t-1$ and $t$ as log differences $r_t=100(\log p_t-\log p_{t-1})$ with $p_t$ closing price on day $t$. In cases of national holidays in some country, the missing index  value is replaced by the last trading  day's value, that is the return is zero. The obtained data sequence has length $N=518$. \begin{figure*}[htbp]
\begin{center}
\includegraphics[width=15cm]{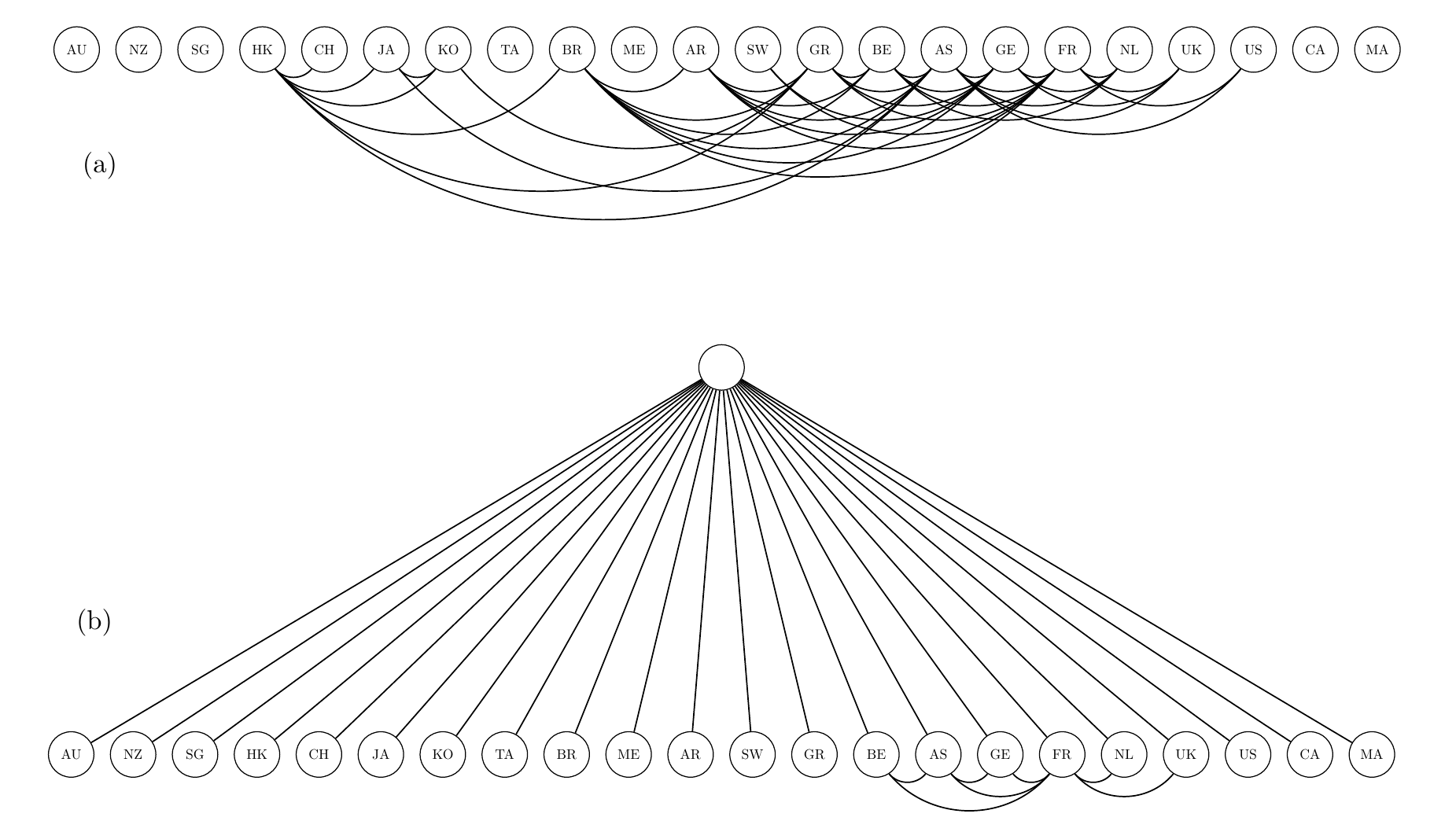}
\end{center}
 \caption{Graphical models for the international financial stock returns data: (a) Best model without latent variables (b) Best  model allowing latent variables.}\label{stock}
\end{figure*}

We applied the identification procedure of Section \ref{susec_synthetic_ex} with $n=1$. In Figure \ref{stock}(b) we depict the estimated graphical model from the financial stock returns data. We found one latent variable and the total number of edges is equal to 29. It is interesting to observe that the latent variable is not sufficient to characterize the conditional dependence relations of Europeans markets (with exception of Greece) and the identification procedure added edges
connecting them. This can be explained by the commencement of the economic and monetary union, see \cite{STOCK_MARKET}. In Figure \ref{stock}(a) we depict the estimated graphical model without latent variables which is characterized by 49 edges among the markets. It is clear that its interpretation is less intuitive than the one with the latent variable. Finally, it is worth noting that $\mathbb{D}(\hat \Phi_C\|\hat\Phi_m^\circ)\cong 3.9$ for both models, that is both models have the same adherence degree to the data.

We consider the estimated joint spectral density $\hat \Phi$ of the manifest and latent variables in (\ref{Phi_p+h_inverse_estimated}) where we choose $\hat\Upsilon_l=1$. Its partial coherence is defined as
\eq \tilde \Phi^{-1}=\mathrm{diag}(\hat\Phi^{-1})^{-1 \slash 2 }\hat\Phi^{-1}\mathrm{diag}(\hat\Phi^{-1})^{-1 \slash 2 }.\eeq
Its entry in position $(k,h)$ represents a measure of how dependent $x_k$ and $x_h$ are conditioned to $\Xc_{V_{m+l}\setminus \{k,h\}}$.
We partition the partial coherence as follows \eq \tilde \Phi^{-1}=\left[
                                                                                              \begin{array}{cc}
                                                                                               \tilde \Upsilon_m & \tilde \Upsilon_{lm}^* \\
                                                                                                \tilde \Upsilon_{lm} & 1 \\
                                                                                              \end{array}
                                                                                            \right].\eeq
In Figure \ref{stock} the entries of $\tilde\Upsilon_{lm}$, representing a measure of the conditional dependence between the latent variable and the stock returns, are depicted.                                                                                              \begin{figure}[htbp]
\begin{center}
\includegraphics[width=9.3cm]{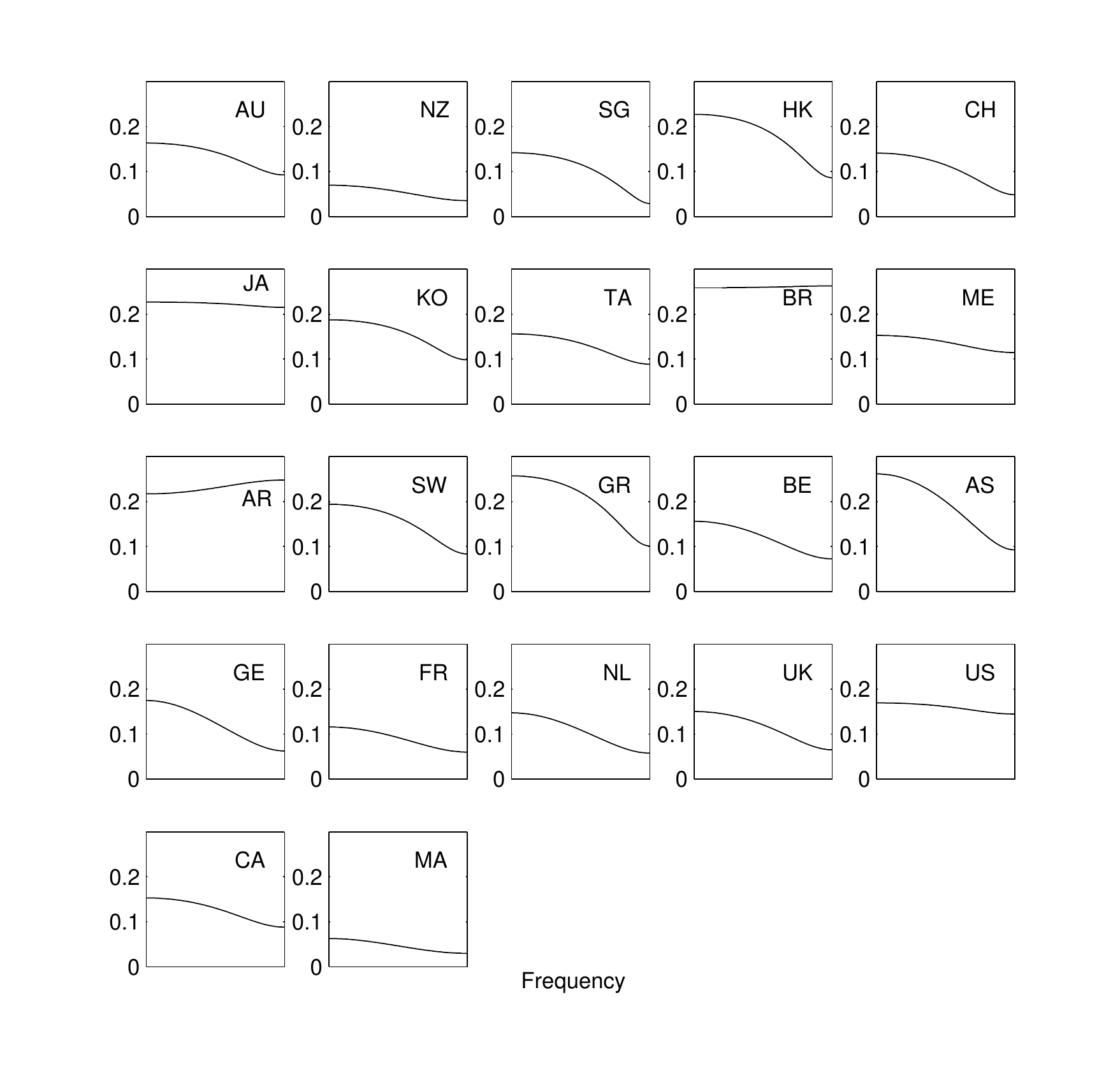}
\end{center}
 \caption{Partial coherence between the latent variable and the stock returns.}\label{partial_coherence}
\end{figure}

\section{Conclusions}
In this paper we dealt with the identification of AR latent-variable graphical models. The inverse of the manifest spectral density of these models
admits a sparse plus low-rank decomposition, captured in two distinct vector subspaces. We presented a two-step procedure for estimating such models.
A first optimization problem uses sparsity and low-rank regularizers to estimate the two vector subspaces. A second optimization problem performs the AR identification
restricted to those vector subspaces. Through duality, the second problem provides a novel covariance extension problem. We provided a simulation study to illustrate the proposed methodology. Finally, we tested our method to international stock return data where the introduction of a latent variable led to a simple graphical model without compromising the adherence to the data.

\appendix

\subsection{Proof of Lemma \ref{lemma_positivity}} \label{proof_lemma_positivity}
If $L\succeq0$, then there exists $C$ such that $L=CC^T$. Accordingly, $\Delta(e^{i\vartheta}) L\Delta(e^{i\vartheta})^*=(\Delta(e^{i\vartheta})C)(\Delta(e^{i\vartheta})C)^*$  which is positive semi-definite for each $\vartheta\in[-\pi,\pi]$. Thus, $\Lambda=\Delta L\Delta^*\succeq 0$. Conversely, if $\Lambda\in\Qc_{m,n}$ is such that $\Lambda\succeq 0$, then it admits the spectral factorization $\Lambda=\Gamma\Gamma^*$ where $\Gamma= \Delta A^T$ such that $A\in\Rs^{m\times m(n+1)}$, \cite{STOICA_MOSES_SPECTRAL_1997}.
Hence, $\Lambda=\Delta A^TA\Delta^*$. We conclude that $\Lambda=\Delta L\Delta^*$ with $L=A^TA\succeq 0$. \qed

\subsection{Proof of Proposition \ref{prop_cov_envelope}}\label{proof_conv_envelope}
Consider an extended-real valued functional $f:\Ac_m\rightarrow [-\infty,+\infty]$. Its {\em conjugate} $f^\star:\Ac_m\rightarrow [-\infty,+\infty]$ is defined as
\eq f^\star (\Phi)=\underset{\Lambda\in\Ac_m}{\sup}( \Sp{\Phi}{\Lambda}-f(\Lambda))\eeq
In view of Theorem 5 in \cite{ROCKA_1974}, the {\em biconjugate} $f^{\star \star}$, i.e. the conjugate of the conjugate, is equal to the convex hull of $f$.

Let $f(\Lambda)=\mathrm{rank}^\prime(\Lambda)$. We prove the statement by showing that $f^{\star \star}$ coincides with (\ref{def_convex_hull}). The proof consists of two steps.

{\em Step 1.} Let $\Dc:=\{\Lambda\in\Ac_m \hbox{ s.t. } \| \Lambda\|\leq 1\}$. Since $f(\Lambda)=+\infty$ for $\Lambda\notin \Dc$, then its conjugate is
\eqn  f^\star(\Phi)&=&\underset{\Lambda \in \Dc}{\sup} (\Sp{\Phi}{\Lambda}-f(\Lambda))\nn \\
&=& \underset{\Lambda \in \Dc}{\sup} \left(\tr\int \Phi \Lambda^*- f(\Lambda)\right)\eeqn where $\Phi\in\Ac_m$.
By applying pointwise the  {\em von Neumann}'s trace theorem \cite{HORN_JOHNSON_MATRIX_AN}, we obtain
\eq \int \tr(\Phi\Lambda^*) \leq \int \sum_{k=1}^m  \sigma_k(\Phi) \sigma_k(\Lambda)\eeq
and equality holds if and only if $\Phi$ and $\Lambda$ admit the following pointwise $SVD$s: $\Phi(e^{i\vartheta})=\Gamma(e^{i\vartheta})\Theta_\Phi(e^{i\vartheta}) \Upsilon(e^{i\vartheta})^*$
and $\Lambda(e^{i\vartheta})=\Gamma(e^{i\vartheta})\Theta_\Lambda(e^{i\vartheta}) \Upsilon(e^{i\vartheta})^*$. Accordingly, $f^\star$ is independent of $\Gamma$ and $\Upsilon$, therefore
\eqn  f^\star(\Phi)&=& \underset{\Lambda\in \Dc}{\sup} \left(  \int \sum_{k=1}^m \sigma_k(\Phi) \sigma_k(\Lambda) - f(\Lambda)\right).\eeqn
If $\Lambda=0$, we have $f^\star(\Phi)=0$ for each $\Phi$. If $f(\Lambda)=l$, with $1\leq l\leq m$, then the supremum is achieved by choosing $\sigma_k(\Lambda(e^{i\vartheta}))=1$ with $k=1\ldots l$, $\vartheta\in[-\pi,\pi]$, and $f^\star(\Phi)=\int \sum_{k=1}^l \sigma_k(\Phi)-l$. Thus, $f^\star$ can be expressed
as \eqn && \hspace{-0.8cm}f^\star(\Phi) =\int \max\left\{ 0, \sigma_1(\Phi(e^{i\vartheta}))-1,\ldots,\sum_{k=1}^l  \sigma_k(\Phi(e^{i\vartheta}))-l,\right.\nn\\  &&\hspace{0.8cm}\left.\ldots,\sum_{k=1}^m  \sigma_k(\Phi(e^{i\vartheta}))-m\right\}\eeqn and the largest term of this set is the one that sums all
positive quantities. We conclude that \eq f^\star(\Phi)= \int \sum_{k=1}^r \left( \sigma_k(\Phi)-1\right),\eeq
where $r(\vartheta)\in \{0,1,\ldots  m\}$ is such that
\eq \label{def_r_theta}\left\{
                         \begin{array}{ll}
                           r(\vartheta)=0, \hbox{ if } \sigma_{1}(\Phi(e^{i\vartheta}))\leq 1 \\
                           \sigma_{r(\vartheta)}(\Phi(e^{i\vartheta}))> 1 \hbox{ and }  \sigma_{r(\vartheta)+1}(\Phi(e^{i\vartheta}))\leq 1,  \hbox{ otherwise.}
                         \end{array}
                       \right.
\eeq In particular, $f^\star(\Phi)=0$ for $\|\Phi\|\leq 1$.

{\em Step 2.} We now compute the conjugate of $f^\star$ which  is defined as
\eq f^{\star\star}(\Lambda)=\underset{\Phi\in\Ac_m}{\sup}( \Sp{\Lambda}{\Phi}-f^\star(\Phi))\eeq where $\Lambda\in\Ac_m$. Proceeding as in Step 1, we have
\eq f^{\star\star}(\Lambda)=\underset{\Phi\in\Ac_m}{\sup}\left( \int \sum_{k=1}^m  \sigma_k(\Lambda)\sigma_k(\Phi)-f^\star(\Phi)\right).\eeq
Next we consider two cases: $\| \Lambda\|> 1$ and $\|\Lambda\|\leq 1$.\\
$\bullet$ {\em Case $\| \Lambda\|>1$.} We have,
\eqn f^{\star\star}(\Lambda)\hspace{-0.3cm}&=&\hspace{-0.3cm}\underset{\Phi\in\Ac_m}{\sup}\left( \int \left(\sum_{k=1}^m \sigma_k(\Lambda)\sigma_k(\Phi)- \sum_{k=1}^r \left( \sigma_k(\Phi)-1\right)\right)\right)\nn\\
&&\hspace{-0.6cm}=  \underset{\Phi\in\Ac_m}{\sup}\left( \int  \left(\sum_{k=1}^r \sigma_k(\Phi)(\sigma_k(\Lambda)-1)\right.\right. \nn\\ &&\left.\left.\hspace{-0.1cm}+ \sum_{k=r+1}^m  \sigma_k(\Phi) \sigma_k(\Lambda)+r\right) \right).\eeqn Let $\bar\vartheta\in [-\pi,\pi] $ such that $\|\Lambda\|=\sigma_1(\Lambda(e^{i\bar\vartheta}))>1$, thus $\sigma_1(\Lambda(e^{i\bar\vartheta}))-1>0$. Since $\Lambda\in\Ac_m$, then $\sigma_k(\Lambda(e^{i\vartheta}))$s are continuous on $\vartheta\in[-\pi,\pi]$, thus  we can choose $\sigma_1(\Phi(e^{i\vartheta}))$ large enough in a neighborhood of $\bar\vartheta$ so that $f^{\star \star}(\Lambda)=+\infty$.\\
$\bullet$ {\em Case $\|\Lambda\|\leq 1$.} If $\|\Phi\|\leq 1$, then $f^\star(\Phi)=0$ and the supremum is achieved by choosing $\Phi=I$, accordingly $\sigma_k(\Phi(e^{i\vartheta}))=1$ for each
$\vartheta\in[-\pi,\pi]$, $k=1\ldots m$, and
\eq f^{\star\star}(\Lambda)=\sum_{k=1}^m \int \sigma_k(\Lambda).\eeq Finally, in the case $\| \Phi\|>1$ the argument of the $\mathrm{sup}$ is always smaller than or equal to the above value:
\eqn && \int \left(\sum_{k=1}^m  \sigma_k(\Lambda)\sigma_k(\Phi)- \sum_{k=1}^r ( \sigma_k(\Phi)-1)\right)\nn\\
 && \hspace{0.5cm} =\int \left(\sum_{k=1}^m  \sigma_k(\Lambda)\sigma_k(\Phi)- \sum_{k=1}^r ( \sigma_k(\Phi)-1)\right.\nn\\ && \hspace{0.9cm} \left. -\sum_{k=1}^m  \sigma_k(\Lambda)\right)+\sum_{k=1}^m \int \sigma_k(\Lambda)\nn\\
 && \hspace{0.5cm} =\int \left(\sum_{k=1}^m  \sigma_k(\Lambda)(\sigma_k(\Phi)-1)- \sum_{k=1}^r ( \sigma_k(\Phi)-1)\right)\nn\\ && \hspace{0.9cm} +\sum_{k=1}^m \int \sigma_k(\Lambda)\nn\\
&& \hspace{0.5cm} =\int \left(\sum_{k=1}^r  \underbrace{(\sigma_k(\Lambda)-1)}_{\hbox{\tiny $\leq 0  \; \vartheta \in [-\pi,\pi]$}}\underbrace{(\sigma_k(\Phi)-1)}_{\tiny \hbox{$>0 \; \vartheta\in[-\pi,\pi]$}}\right.\nn\\
&&\hspace{0.9cm}+ \left. \sum_{k=r+1}^m \sigma_k(\Lambda)\underbrace{( \sigma_k(\Phi)-1)}_{\hbox{\tiny $\leq 0 \;\vartheta\in[-\pi,\pi]$}}\right) +\sum_{k=1}^m \int \sigma_k(\Lambda)\nn\\ &&  \hspace{0.5cm} \leq  \sum_{k=1}^m \int \sigma_k(\Lambda)\eeqn
where we exploited (\ref{def_r_theta}).

We conclude that \eq f^{\star\star}(\Lambda)=\left\{
  \begin{array}{ll}
    \sum_{k=1}^m \int \sigma_k(\Lambda), & \|\Lambda \|\leq 1 \\
    +\infty, & \hbox{otherwise.}
  \end{array}
\right.\eeq \qed

\subsection{Proof of Proposition \ref{prop_dual_dual_Delta_regularized}}\label{proof_prop_dual_dual_Delta_regularized}
Before proving the statement, we establish the following lemma.
\lem \label{lemma_P_Omega} Let $\mathbf{C}$ be a closed convex subset of $\{Z\in \Mb_{m,n} \hbox{ s.t. } \tr(Z_0)=0\}$,
$c$ be a constant term. If the following convex optimization problem
is feasible \eqn &\underset{\substack{W\in\Qb_m\\ Z\in\Mb_{m,n}}}{\max} &\log\det W +c\nn\\
&\hbox{subject to}& W\succ 0\nn\\
&& \Td(\hat R)+\mathrm{T}(Z)\succeq \left[
                                 \begin{array}{cc}
                                   W & 0 \\
                                   0 & 0 \\
                                 \end{array}
                               \right]
\nn \nn\\
&& Z\in \mathbf{C}  \eeqn
then it admits a solution.\elem \IEEEproof  By assumption, the optimization problem is feasible, i.e. there exist $\bar W\in\Qb_m$ and $\bar Z\in\Mb_{m,n}$ satisfying the constraints, and such that $|\log\det \bar W +c|<\infty$. Accordingly, the above problem is equivalent to maximize $\log\det W$ over the set \eqn && \hspace{-0.7cm}\mathbf{D}:=\left\{(W,Z)\in \Qb_m\times \mathbf{C} \hbox{ s.t. }  W\succ 0,\right.\nn\\ && \left.  \Td(\hat R)+\mathrm{T}(Z)\succeq \left[
                                 \begin{array}{cc}
                                   W & 0 \\
                                   0 & 0 \\
                                 \end{array}
                               \right],\;\;  \log\det W\geq \log\det \bar W    \right\}.\nn\eeqn
Next we show that $\mathbf{D}$ is a compact set. Since $\log\det W $ is continuous over $\mathbf{D}$, it follows from {\em Weierstrass'} theorem that $\log\det W$ admits a maximum on $\mathbf{D}$.

To prove the compactness of $\mathbf{D}$, we show that it is bounded and closed. Let $\{(Z^{(k)},W^{(k)})\}_{k\in\Ns}$
be a sequence belonging to $\mathbf{D}$. Since the minimum singular value of the map $\Td$
is strictly positive, if $\| Z^{(k)}\|\rightarrow \infty$ as $k\rightarrow \infty$, then $\|T(Z^{(k)})\|\rightarrow +\infty$. Since $\mathrm{T}(Z^{(k)})$ is a symmetric matrix, $\mathrm{T}(Z^{(k)})$ has at least one eigenvalue tending to infinity in modulus. Moreover $\tr(\mathrm{T}(Z^{(k)}))=0$ because $Z\in \mathbf{C}$. Thus
 $\mathrm{T}(Z^{(k)})$, and hence $\Td(\hat R) +\mathrm{T}(Z^{(k)})$, has at least one eigenvalue tending to $-\infty$. This is not possible because $Z^{(k)}$ must satisfy inequality \eq \Td(\hat R)+\mathrm{T}(Z^{(k)}) \succeq \left[
                                                                                                         \begin{array}{cc}
                                                                                                           W & 0 \\
                                                                                                           0 & 0\\
                                                                                                         \end{array}
                                                                                                       \right]
\succeq 0.\eeq Thus, $\|Z^{(k)}\|<\infty$. Moreover, $\| W^{(k)}\|<\infty$ because $0\prec W^{(k)}\preceq ( \Td(\hat R)+\mathrm{T}(Z^{(k)}))_{00}$. Therefore $\mathbf{D}$ is bounded. Let $\partial \mathbf{D} $ denote the subset of the boundary of $\mathbf{D}$ not contained in $\mathbf{D}$. Since $\mathbf{C}$ is a closed subset of $\Mb_{m,n}$, $\partial \mathbf{D}$ is at most the set of elements $(Z,W)$ such that $W$ is positive semi-definite and singular. Since $\lim _{(Z,W)\rightarrow \partial \mathbf{D}} \log \det W=-\infty$ and $W$ must satisfy the inequality $\log\det W \geq \log\det \bar W $, we conclude that $\partial \mathbf{D}$ is an empty set. Accordingly, $\mathbf{D}$ is closed. \qed\\

We proceed to prove Proposition \ref{prop_dual_dual_Delta_regularized}. Since $\Td(\hat R)\succ0$,
Problem (\ref{PB:dual_dual_Delta_regularized}) is feasible (it is sufficient to pick $W=\alpha I$ and $Z=0$ where $\alpha>0$ is the minimum eigenvalue of $\Td(\hat R)$). Then, by applying Lemma \ref{lemma_P_Omega} with
\eqn &&\hspace{-0.6cm}\mathbf{C}:= \{Z\in\Mb_{m,n} \hbox{ s.t. } \mathrm{diag}(Z_j)=0\;j=0\ldots n,\nn\\ &&  \hspace{0.1cm}\sum_{j=0}^n |(Z_j)_{kh}|+|(Z_j)_{hk}|\leq \lambda\gamma\; k\neq h,\; \lambda I+\mathrm{T}(Z)\succeq 0  \}\nn\eeqn
we conclude that (\ref{PB:dual_dual_Delta_regularized}) admits a solution. Finally, it is worth noting the objective function in (\ref{PB:dual_dual_Delta_regularized}) is strictly convex with respect to $W$, thus the optimal solution $W^\circ$ is unique. \qed

\subsection{Proof of Proposition \ref{prop_equivalence}}\label{proof_prop_equivalence}

Our proof uses the following lemma whose proof can be found in  \cite{SONGSIRI_GRAPH_MODEL_2010}.
\lem \label{lemma_basics_songsiri} Let $Z\in\Mb_{m,n}$, $W\in\Qb_m$. If $W\succ 0$ and such that \eq \Td(Z)\succeq \left[
                                                                                            \begin{array}{cc}
                                                                                              W & 0 \\
                                                                                              0 & 0 \\
                                                                                            \end{array}
                                                                                          \right]
\eeq then $\Td(Z)\succ 0$ and the unique solution to the {\em Yule-Walker equations}, \cite{SODDERSTROM_STOICA_SYS_ID}, \eq \left\{
                                                                                              \begin{array}{ll}
                                                                                                \Td(Z)B^T=\left[
                                                                                                                 \begin{array}{c}
                                                                                                                   W \\
                                                                                                                   0 \\
                                                                                                                 \end{array}
                                                                                                               \right],& B\in \Rs^{m \times m(n+1)} \\
                B_0=I                                                                                              \end{array}
                                                                                            \right.
\eeq  is such that $B \Delta^*$ has zeros inside the unit circle.\elem

We proceed to prove Proposition \ref{prop_equivalence}. Note that the duality gap between (\ref{PB:dual_Delta_regularized}) and (\ref{PB:dual_dual_Delta_regularized}) is equal to zero. In particular, \eq \label{trace_Qo_DeltaO} \Sp{U^\circ}{X^\circ}=0\eeq
where $U^\circ\in\Qb_{m(n+1)}$, $U^\circ\succeq 0$ maximizes (\ref{funzionale_duale_regolarizzato}). Note that $U^\circ$ can be expressed in the following way
\eq U^\circ=\Td(\hat R)+\Td(Z^\circ)-\left[
                                                        \begin{array}{cc}
                                                          W^\circ & 0 \\
                                                          0 & 0 \\
                                                        \end{array}
                                                      \right]
\eeq
where $W^\circ\succ 0$ and $Z^\circ\in\Mb_{m,n}$ are solution to Problem (\ref{PB:dual_dual_Delta_regularized}). By Lemma \ref{lemma_basics_songsiri}, we have that  $\Td(\hat R)+\Td(Z^\circ)\succ 0$, accordingly $U^\circ$ has rank at least equal to $mn$. Since $U^\circ,X^\circ\succeq 0$, (\ref{trace_Qo_DeltaO}) implies that $X^\circ$ has rank at most equal to $m$. On the other hand $\mathrm{rank}(X^\circ)\geq m$ because $X_{00}^\circ=(W^\circ)^{-1}\succ0$. We conclude that $\mathrm{rank}(X^\circ)=m$. Hence, there exists $A\in\Rs^{m \times m(n+1)}$ full row rank
such that $X^\circ=A^TA$ with $X_{00}^\circ=A_0^T A_0$. Since $U^\circ,X^\circ\succeq 0$, (\ref{trace_Qo_DeltaO}) implies
\eq \left(\Td(\hat R)+\Td(Z^\circ)-\left[
                                               \begin{array}{cc}
                                                 W^\circ & 0 \\
                                                 0 & 0 \\
                                               \end{array}
                                             \right]\right)A^T=0.\eeq
By defining $B\in\Rs^{m \times m(n+1)}$ such that $B=A_0^{-1}A$ we obtain
\eq\label{Yule_walker_sharp} (\Td(\hat R)+\Td(Z^\circ))B^T=\left[
                                               \begin{array}{c}
                                                 W^\circ\\
                                                 0  \\
                                               \end{array}
                                            \right],\;\; B_0=I .\eeq
Since  $\Td(\hat R)+\Td(Z^\circ)\succ 0 $, the Yule-Walker equations (\ref{Yule_walker_sharp}) admits a unique solution such that $B \Delta^*$ has zeros inside the unit circle. Accordingly, there exists $X^\circ$ such that
\eq \Delta X^\circ \Delta^*=\Delta A^T A \Delta^*=(\Delta B^T) (W^\circ)^{-1}(B \Delta^*)\succ 0.\eeq
Finally, uniqueness of $ X^\circ$ follows from the uniqueness of $W^\circ$ and $B$. It remains to be shown the existence of $L^\circ$. In view of (\ref{PB:dual_Delta_regularized}), we have \eqn L^\circ=&\hspace{-0.5cm}\underset{\substack{L\in\Qb_{m(n+1)}}}{\arg\min}& \hspace{-0.5cm}\lambda\gamma h_\infty(\Dd(X^\circ+L)) +\lambda\tr(L)\nn\\
   & \hspace{-0.5cm}\hbox{ subject to }&    \hspace{-0.5cm} \;L\succeq 0\eeqn
where the objective function is continuous. Since $L=0$ is a feasible point, we can restrict $L$ to belong to
\eqn &&\hspace{-0.8cm}\mathbf{D}:=\{L\in\mathbf{Q}_{m(n+1)} \hbox{ s.t. } L\succeq 0, \nn\\
&& \hspace{-0.6cm}\lambda\gamma h_\infty(\Dd(X^\circ+L)) +\lambda\tr(L)\leq \lambda\gamma h_\infty(\Dd(X^\circ))\}.\eeqn
It is not difficult to show that $\mathbf{D}$ is a closed and bounded set, therefore by {\em Weierstrass'} theorem $L^\circ$ does exist. \qed\\

\subsection{Proof of Proposition \ref{prop_transversality}} \label{proof_prop_transversality}
By Proposition \ref{prop_equivalence}, $X^\circ$ is unique. If (\ref{constr_Sigma}) admits a unique solution $H$, then $L^\circ =G^T H G $ is unique. Therefore,
$\Vc_{E_m}$ and $\Vc_G$ are unique because the uniqueness of $X^\circ$ and $L^\circ$. Equation (\ref{constr_Sigma}) may be written in the compact form \eq Ay=b\eeq
where the vector $y\in\Rs^{l(l+1)/2}$ contains the independent parameters of $H$, $A\in\Rs^{(n+1)|I|\times l(l+1)/2}$ only depends on $G$ and $b\in\Rs^{(n+1)|I|}$ only depends on $X^\circ$. If (\ref{constr_Sigma}) admits a unique solution, then it is obtained in the following way \eq y=(A^T A)^{-1}A^T b\eeq
and changing $b$ (i.e. $X^\circ$) such a solution is still unique. Accordingly, the uniqueness of the solution to (\ref{constr_Sigma})
is equivalent to the uniqueness of the decomposition \eq \Phi^{-1}_m=\Sigma-\Lambda\eeq with $\Phi_m^{-1}\in \Vc_{E_m}+\Vc_G$, $\Sigma\in\Vc_{E_m}$ and $\Lambda\in\Vc_G$. Therefore, $\Vc_{E_m} \cap \Vc_G=\{0\} $.\qed

\subsection{Proof of Proposition \ref{prop_existence_sol}} \label{proof_existence_sol}
The statement can be proved by duality theory along the same line of the proof of Proposition \ref{prop_dual_dual_Delta_regularized} and Proposition \ref{prop_equivalence}, respectively.\qed

\end{document}